\documentclass[reqno,twoside,11pt]{amsart}
 
\hoffset - 1.5cm

\usepackage{graphicx}
\usepackage{amstext}
\usepackage{pstricks,pst-node,pst-text,pst-3d}
\usepackage{color}
\usepackage{amsmath}
\usepackage{amsthm}
\usepackage{amssymb}
\newtheorem{Theorem}{Theorem}[section]
\newtheorem{Fact}{Fact}[section]
\newtheorem{Lemma}{Lemma}[section]
\newtheorem{Proposition}{Proposition}[section]
\newtheorem{Corollary}{Corollary}[section]
\theoremstyle{definition}
\newtheorem{Definition}{Definition}[section]
\newtheorem{Example}{Example}[section]
\newtheorem{Remark}{Remark}[section]

\newcommand{\ba}{\begin{array}}
\newcommand{\bc}{\begin{center}}
\newcommand{\bd}{\begin{description}}
\newcommand{\bdm}{\begin{displaymath}}
\newcommand{\be}{\begin{enumerate}}
\newcommand{\beq}{\begin{equation}}
\newcommand{\bdf}{\begin{Definition}}
\newcommand{\bex}{\begin{Example}}
\newcommand{\bft}{\begin{Fact}}
\newcommand{\bl}{\begin{Lemma}}
\newcommand{\bp}{\begin{Proposition}}
\newcommand{\br}{\begin{Remark}}
\newcommand{\bt}{\begin{Theorem}}
\newcommand{\bco}{\begin{Corollary}}
\newcommand{\bh}{\begin{Hipothesis}}
\newcommand{\ea}{\end{array}}
\newcommand{\ec}{\end{center}}
\newcommand{\ed}{\end{description}}
\newcommand{\edm}{\end{displaymath}}
\newcommand{\ee}{\end{enumerate}}
\newcommand{\eeq}{\end{equation}}
\newcommand{\edf}{\end{Definition}}
\newcommand{\eex}{\end{Example}}
\newcommand{\eft}{\end{Fact}}
\newcommand{\el}{\end{Lemma}}
\newcommand{\ep}{\end{Proposition}}
\newcommand{\er}{\end{Remark}}
\newcommand{\et}{\end{Theorem}}
\newcommand{\eco}{\end{Corollary}}
\newcommand{\eh}{\end{Hipothesis}}

\newcommand{\bA}{\mathbb{A}}

\newcommand{\bH}{\mathbb{H}}
\newcommand{\bI}{\mathbb{I}}

\newcommand{\bN}{\mathbb{N}}

\newcommand{\bR}{\mathbb{R}}

\newcommand{\bT}{\mathbb{T}}

\newcommand{\bV}{\mathbb{V}}
\newcommand{\bW}{\mathbb{W}}
\newcommand{\bX}{\mathbb{X}}
\newcommand{\bY}{\mathbb{Y}}
\newcommand{\bZ}{\mathbb{Z}}
\newcommand{\cC}{\mathcal{C}}

\newcommand{\cE}{\mathcal{E}}
\newcommand{\cF}{\mathcal{F}}

\newcommand{\cI}{\mathcal{I}}

\newcommand{\cK}{\mathcal{K}}
\newcommand{\cL}{\mathcal{L}}

\newcommand{\cN}{\mathcal{N}}
\newcommand{\cO}{\mathcal{O}}

\newcommand{\cT}{\mathcal{T}}
\newcommand{\cU}{\mathcal{U}}

\newcommand{\numsec}{\setcounter{Theorem}{0}\setcounter{Definition}{0}
\setcounter{Remark}{0} \setcounter{Lemma}{0} \setcounter{Fact}{0}
\setcounter{Proposition}{0} \setcounter{Corollary}{0}
\setcounter{Example}{0} \setcounter{equation}{0}
\setcounter{Property}{0}\renewcommand\theequation{\arabic{section}.\arabic{equation}}
\renewcommand\theTheorem{\arabic{section}.\arabic{Theorem}}
\renewcommand\theDefinition{\arabic{section}.\arabic{Definition}}
\renewcommand\theRemark{\arabic{section}.\arabic{Remark}}
\renewcommand\theLemma{\arabic{section}.\arabic{Lemma}}
\renewcommand\theFact{\arabic{section}.\arabic{Fact}}
\renewcommand\theProposition{\arabic{section}.\arabic{Proposition}}
\renewcommand\theCorollary{\arabic{section}.\arabic{Corollary}}
\renewcommand\theExample{\arabic{section}.\arabic{Example}}
\renewcommand\theProperty{\arabic{section}.\arabic{Property}}}

\numberwithin{equation}{section} \errorcontextlines=0

\newcommand{\sone}{S^1}
\newcommand{\sotwo}{SO(2)}
\newcommand{\otwo}{O(2)}
\newcommand{\sot}{SO(3)}
\newcommand{\ds}{\displaystyle}
\newcommand{\nt}{\noindent}

\newcommand{\h}{\mathbb{H}}

\newcommand{\on}{O(n,\bR)}

\newcommand{\vp}{\varphi}

\newcommand{\lin}{\mathrm{span}}
\newcommand{\subg}{\overline{\mathrm{sub}}(G)}
\newcommand{\subgc}{\overline{\mathrm{sub}}[G]}
\newcommand{\subh}{\overline{\mathrm{sub}}(H)}
\newcommand{\subhc}{\overline{\mathrm{sub}}[H]}
\newcommand{\subt}{\overline{\mathrm{sub}}(\bT)}

\newcommand{\subotwo}{\overline{\mathrm{sub}}(\otwo)}

\newcommand{\sub}{\overline{\mathrm{sub}}}
\newcommand{\card}{\mathrm{card}\:}
\newcommand{\morse}{\mathrm{m^-}}

\newcommand{\chig}{\chi_G}
\newcommand{\chih}{\chi_H}
\newcommand{\chihn}{\chi_{H^{\nu}}}

\newcommand{\ci}{\cC\cI}
\newcommand{\cig}{\ci_G}
\newcommand{\cih}{\ci_H}
\newcommand{\cihn}{\ci_{H^{\nu}}}

\usepackage{color}

\begin{document}

\title[Periodic orbits]{Symmetric Liapunov center theorem}

\author{Ernesto P\'{e}rez-Chavela$^{1)}$}
\address{$^{1)}$ Departamento de Matemáticas \\ Instituto Tecnológico Autónomo de M\'{e}xico (ITAM), Río Hondo 1, Col. Progreso Tizap\'{a}n
M\'{e}xico D.F. 01080 \\ M\'{e}xico}

\author{S{\l}awomir Rybicki$^{2)}$}
\address{$^{2)}$ Faculty of Mathematics and Computer Science\\
Nicolaus Copernicus University \\
PL-87-100 Toru\'{n} \\ ul. Chopina $12 \slash 18$ \\
Poland}

\author{Daniel Strzelecki$^{2)}$}

\email{ernesto.perez@itam.mx (E. P\'{e}rez-Chavela)}
\email{rybicki@mat.umk.pl (S. Rybicki)}
\email{danio@mat.umk.pl (D. Strzelecki)}

\date{\today}

\thanks{$^{1)}$ Partially supported by the Asociación Mexicana de Cultura A.C}
\thanks{$^{2)}$ Partially supported   by the National Science Center,  Poland,  under grant    DEC-2012/05/B/ST1/02165.}

\begin{abstract}

In this article, using an infinite-dimensional equivariant Conley index, we prove a~generalization of the profitable Liapunov center theorem for symmetric potentials. Consider the system $\ddot{q}= -\nabla U(q),$ where $U(q)$ is a
$\Gamma$-symmetric potential, where $\Gamma$ is a compact Lie group  acting linearly on $\bR^n$. If the system possess a non-degenerate orbit  of stationary solutions  $\Gamma(q_0)$ with trivial isotropy group , such that there exists at least one positive eigenvalue of the Hessian $\nabla^2 U(q_0)$,  then in any neighborhood of  orbit $\Gamma(q_0)$ there is a periodic orbit of solutions of  the system.
\end{abstract}

\maketitle

%%%%%%%%%%%%%%%%%%%%%%%%%%%%%%%%%%%%%%%%%%%%%%%%%%%%%%%%%%%%%%
\numsec
\section{Introduction}

One of the most famous theorems concerning the existence of periodic solutions of ordinary differential equations is the celebrated Liapunov center theorem \cite{[H]}. Consider the second order autonomous system $\ddot q(t)  =  - \nabla U(q(t)),$ where $U \in C^2(\bR^n,\bR), \nabla U(0)=0$ and $ \det \nabla^2 U(0) \neq 0.$ Let be $\sigma(\nabla^2 U(0))$ the spectrum of the respective Hessian at $0$, 
Liapunov's center theorem says that if  $\sigma(\nabla^2 U(0)) \cap (0,+\infty) = \{\beta_1^2,\ldots, \beta_m^2\}$ for $\beta_1 >\ldots > \beta_m > 0$ and there is $\beta_{j_0}$ satisfying $\beta_1 \slash \beta_{j_0}, \ldots, \beta_{j _0-1} \slash \beta_{j_0} \not \in \bN,$ then  there is a sequence $\{q_k(t)\}$ of periodic solutions of system  
\beq \label{newsys} \ddot q(t)  =  - \nabla U(q(t)),\eeq 
with amplitude  tending to zero and the  minimal period tending to $2 \pi \slash \beta_{j_0}.$ Proof of this theorem one can find in \cite{[MAWI]}, see also \cite{[BERGER1],[BERGER2],[SZULKIN]}.

Generalization of this theorem in two directions is due to Szulkin \cite{[SZULKIN]}.  In the first direction Szulkin, using the infinite-dimensional Morse theory for strongly-indefinite functionals, proved Liapunov type center theorem for Hamiltonian systems   
\beq \label{hamsys} \dot z(t)=J \nabla H(z(t)),\eeq 
where $H \in C^2(\bR^{2n},\bR), \nabla H(0)=0,\det \nabla^2 H(0) \neq 0,$ and $J=\left[\ba{rr} 0 & -I \\ I & 0 \ea \right]$ is the standard simplectic matrix in $\bR^{2n}$, see Theorem 4.1 and Corollaries 4.2, 4.3 of \cite{[SZULKIN]}. The next important generalization  was to consider system \eqref{hamsys} with Hamiltonian $H$ for which $0 \in \bR^{2n}$ is an isolated degenerate critical point of the Hamiltonian $H$ with nontrivial Conley index i.e. $\nabla H(0)=0,$ $ 0\in \bR^{2n}$ is isolated in  $(\nabla H)^{-1}(0)$ and $\ci(\{0\},-\nabla H) \neq [*,*],$  where $\ci$ denotes the Conley index and $*$ denotes a point. The Liapunov center type theorem for such system has also been proved in \cite{[SZULKIN]}, see theorem 4.4 and corollary 4.5 of \cite{[SZULKIN]}. 

Taking the advantage of  variational structure of system \eqref{hamsys}, usually the authors convert the problem of the existence of non-stationary periodic solutions in a neighborhood of a stationary one into a bifurcation problem. Finally, they  apply the Morse theory or the Conley index theory to prove the existence of bifurcation of periodic solutions. In this way they obtain a local bifurcation (a sequence of solutions bifurcating from the family of trivial ones) which does not have to be global (a connected set of solutions bifurcating from the family of trivial ones), see \cite{[AMB],[BOHME],[IZE], [MARINO], [TAKENS]} for discussions and examples.

Another generalization  of the Liapunov center theorem is due to Dancer and the second author  \cite{[DARY]}. They considered system \eqref{hamsys}  with Hamiltonian $H$ for which $0 \in \bR^{2n}$ is an isolated degenerate critical point of $H$ with nontrivial Brouwer  index i.e.  $0 \in \bR^{2n}$ is isolated in $(\nabla H)^{-1}(0)$ and $\deg_{B} (\nabla H,B^{2n}_{\alpha},0) \neq  0,$ where $\alpha > 0$ is sufficiently small and $\deg_{B}(\cdot)$ is the Brouwer degree. Note that since  $\chi(\ci(\{0\},-\nabla H)) = \deg_{B} (\nabla H,B^n_{\alpha},0),$ the assumption considered in \cite{[DARY]} implies that of \cite{[SZULKIN]}, where $\chi(\cdot)$ is the Euler characteristic. Under this stronger assumption they have proved that there is a connected set of non-stationary periodic solutions of system \eqref{hamsys} emanating from the stationary solution $u_0 \equiv 0.$ In order to prove this theorem
they have applied the degree theory for $\sone$-equivariant gradient maps, see \cite{[GEBA]}.

There are also theorems giving estimations of the number of periodic orbits of system \eqref{hamsys} on an energy level close to the non-degenerate critical point $0 \in \bR^{2n}$of the Hamiltonian $H$ due to Weinstein \cite{[WEINSTEIN]}  and Moser \cite{[MOSER]}. For differential equations with first integral there are similar results due to Dancer and Toland \cite{[DATO]}, Marzantowicz and Parusi\'{n}ski \cite{[MARPAR]}.

It can happen that the stationary solutions of system \eqref{newsys} are not isolated critical points of the potential $U$ and the set of stationary solutions consists of the orbits of a compact Lie group $\Gamma.$
For example the Lennard-Jones potential $U : \Omega \to \bR$ is $\Gamma=\sotwo$-invariant and $(\nabla U)^{-1}(0) \cap \Omega$ consist of $\Gamma$-orbits i.e. the stationary solutions of system \eqref{newsys} are not isolated, see \cite{[CLC1],[CLC2]}. It is worth pointing out that  one can not apply the theorems mentioned above to the study of  non-stationary periodic solutions of system \eqref{newsys}. 

The goal of  this paper is to prove  a symmetric version of the Liapunov center theorem. Let  $\Omega \subset \bR^n$ be an open and $\Gamma$-invariant subset of an orthogonal   representation $\bR^n$ of a compact Lie group $\Gamma.$   Assume that $q_0 \in \Omega$ is a critical point of  a~$\Gamma$-invariant potential $U : \Omega \to \bR$  of  class $C^2$ with an isotropy group  $\Gamma_{q_0}=\{\gamma \in  \Gamma : \gamma q_0=q_0\}.$  Since the gradient $\nabla U : \Omega \to \bR^n$ is $\Gamma$-equivariant, the $\Gamma$-orbit $\Gamma(q_0)=\{\gamma q_0 : \gamma \in \Gamma\}$ consists of critical points of the potential $U$ i.e. $\Gamma(q_0) \subset (\nabla U)^{-1}(0),$ and therefore  $\dim \ker \nabla^2 U(q_0) \geq \dim \Gamma(q_0).$

The main result of this article is the following.
\bt\label{main-theo} [Symmetric Liapunov center theorem]  Under the above assumptions. If moreover,
\begin{enumerate}
\item  the isotropy group $\Gamma_{q_0}$ is trivial,
\item  $\dim \ker \nabla^2 U(q_0) = \dim \Gamma(q_0),$
\item there exists at least one positive eigenvalue of the Hessian $\nabla^2 U(q_0).$
\end{enumerate}
Then in any open neighborhood of the orbit $\Gamma(q_0)$ there is an orbit of non-stationary periodic solutions of system $\ddot q(t)=-\nabla U(q(t))$. Moreover,   if   $\sigma(\nabla^2 U(q_0)) \cap (0,+\infty) = \{\beta_1^2,\ldots, \beta_m^2\}$ with $\beta_1 >\ldots > \beta_m > 0$ and there exists a fixed $\beta_{j_0}$ satisfying $\beta_1 \slash \beta_{j_0}, \ldots, \beta_{j _0-1} \slash \beta_{j_0} \not \in \bN,$ then the  minimal periods of periodic orbits in the neighborhood of the orbit $\Gamma(q_0)$ are close to $2 \pi \slash \beta_{j_0}.$ 
\et
To prove this theorem we apply the infinite-dimensional version of the $(\Gamma \times \sone)$-equivariant Conley index theory due to Izydorek \cite{[IZYDOREK]}. We emphasize that if the group $\Gamma$ is trivial then the above theorem is the classical Liapunov center theorem.

After this introduction our article is organized as follows. 
In the first part of Section \ref{prelim} we consider $2 \pi$-periodic solutions of system \eqref{newsys} as critical points of a  functional $\Phi$ defined on a~suitably chosen Hilbert space $\bH^1_{2\pi}$. Additionally, we present some  properties of the Hessian of this functional. Next we summarize without proofs the relevant material on equivariant topology and representation theory of compact Lie groups. We have introduced  the notion of an admissible $(G,H)$ pair of compact Lie groups, where  $H \in \subg,$ see Definition  \ref{admp} and the Euler ring $U(G)$ of a compact Lie group $G,$ see Definition \ref{tdr} and Lemma \ref{ms}. The special case of the Euler ring $U(\sone)$ is discussed in Remark \ref{invert}.  A formula for the $G$-equivariant Euler characteristic $\chig(X) \in U(G)$ of a finite pointed $G$-CW-complex $X$ is presented in Lemma \ref{pc}.  
In Theorems \ref{hsm}, \ref{hsmadm} we have expressed the $G$-equivariant Euler characteristic $\chig(G^+ \wedge_H X) \in U(G),$ of a~$G$-CW-complex $G^+ \wedge_H X$ in therms of the $H$-equivariant Euler characteristic $\chih(X) \in U(H)$ of a $H$-CW-complex $X.$ We underline that if the pair $(G,H)$ is admissible then the map $U(H) \ni \chih([X]_H) \to \chig\left([G^+ \wedge_H X]_G\right) \in  U(G)$ is injective, see Theorem \ref{hsmadm} and Corollary \ref{hsmadmco}.

Section \ref{conley}   is devoted to the computations of the $G$-equivariant Conley index $\cig(G(x_0),-\nabla \vp)$ of a  non-degenerate orbit $G(q_0)$ of critical points of an invariant potential $\vp \in C^2_G(\Omega,\bR).$  First of all we are interested in finding relation between the equivariant Conley index of a non-degenerate orbit and the equivariant Conley index of a non-degenerate critical point of the potential restricted to the space orthogonal to this orbit. Such relation is proved in Theorem \ref{cio}. This relation allows us to distinguish the equivariant Conley index  of non-degenerate orbits analyzing only the potentials restricted to the  orthogonal spaces to these orbits, see Corollaries \ref{cio1},  \ref{hgg}, \ref{hggr}.  Finally in  Theorem \ref{last} we distinguish  equivariant Conley index of so called special non-degenerate orbits.  

 Our main results are proved in Section \ref{main},  where we consider system \eqref{newsys} with $\Gamma$-symmetric potential $U$ and   study periodic solutions  of this system in a neighborhood of a non-degenerate orbit $\Gamma(q_0)$ of critical points of $U$ i.e. we have proved the Symmetric  Liapunov center theorem (Theorem \ref{main-theo}). This Theorem is a natural generalization of the classical Liapunov center theorem. The basic idea is to consider periodic solutions of system \eqref{newsys}  as critical orbits of $G=(\Gamma \times \sone)$-invariant family of functionals, see equation  \eqref{gradf}. In other words we have converted the problem of the existence of periodic solutions of system \eqref{newsys} in a neighborhood of $\Gamma(q_0)$ into $G$-symmetric, infinite-dimensional  and variational bifurcation problem. 

To prove the existence of bifurcation we use the infinite-dimensional equivariant Conley index due to Izydorek \cite{[IZYDOREK]}. First we prove a technical lemma \ref{cis} which yields information on the $\sone$-equivariant Conley indices of  critical points of functionals restricted to the  orthogonal space to the orbit $\Gamma(q_0) \subset \bH^1_{2\pi}.$ Next we prove Theorem \ref{main-theo}. 

Finaly in Section \ref{applications} we consider two simple examples coming from celestial mechanics just to show  the strength of our main result, and how we can use it. Specifically we analyze a couple of generic galactic type potentials, and show how to find periodic orbits on them.

%%%%%%%%%%%%%%%%%%%%%%%%%%%%%%%%%%%%%%%%%%%%%%%%%%%%%%%%%%%%%%
\numsec
\section{Preliminaries}
\label{prelim}

In this section we give a brief exposition of material on functional analysis which we will need in the rest of this article. Our purpose is summarize without proofs the relevant tools on  equivariant topology used along this paper.

Fix  an open set $\Omega \subset \bR^n$ and consider the following system of second order equations

\beq \label{sys}  
\left\{
\ba{rcl}  \ddot q(t) & = & - \nabla U(q(t)) \\
q(0) & = & q(2\pi) \\
\dot q(0) & = & \dot q(2\pi)
\ea,
\right.
\eeq 
where $U \in C^2(\Omega, \bR).$

We define 
$$\h^1_{2\pi} = \{u : [0,2\pi] \rightarrow \bR^n : \text{ u is abs. continuous map, } u(0)=u(2\pi), \dot u \in L^2([0,2\pi],\bR^n)\}$$ and an open subset $\h^1_{2\pi}(\Omega) \subset \h^1_{2\pi} $ by 
$\ds \h^1_{2\pi}(\Omega)=\{ u \in \h^1_{2\pi}  : u([0,2\pi]) \subset \Omega\}.$

 It is well known that $\h^1_{2\pi}$ is a separable Hilbert space with a scalar product given by the
formula $$\ds \langle u,v\rangle_{\h^1_{2\pi}} = \int_0^{2\pi} (\dot u(t), \dot v(t)) + (u(t),v(t)) \; dt,$$ where $(\cdot,
\cdot)$ and $\| \cdot \|$ are the usual scalar product and norm in  $\bR^n,$ respectively. It is easy to show
that $\left(\h^1_{2\pi},\langle\cdot,\cdot\rangle_{\h^1_{2\pi}}\right)$ is an orthogonal representation of the group $\sone$ with an $\sone$-action given by shift in time. It is clear that $\h^1_{2\pi}(\Omega)$ is $\sone$-invariant.

Let be $ \{ e_1,\ldots,e_n \} \subset \bR^n$ be the standard basis in $\bR^n.$ Define $\h_0=\bR^n, \h_k= \lin \{e_i \cos kt, e_i \sin kt: i=1,\ldots,n\}$ and note that 
\beq \label{space} 
\bH^1_{2\pi} = \overline{\h_0 \oplus \bigoplus_{k=1}^{\infty} \h_k} 
\eeq
and that the finite-dimensional spaces $\h_k,k=0,1,\ldots$ are orthogonal representations of $  \sone.$
Define an $S^1$-invariant functional 
$\Phi : \h^1_{2\pi}(\Omega)  \to \bR$ of the class $C^2$ as follows
$$
\Phi(q) = \int_0^{2\pi} \left( \frac{1}{2} \| \dot q(t) \|^2 -  U(q(t)) \right) \; dt,
$$
notice that for any $q\in \h^1_{2\pi}(\Omega)$ and  $q_1 \in \h^1_{2\pi}$ we have $D\Phi(q)(q_1)=\langle \nabla \Phi(q),q_1\rangle_{\h^1_{2\pi}}=
\langle q - \nabla \zeta(q),q_1\rangle_{\h^1_{2\pi}},$ where $\nabla \zeta : \h^1_{2\pi}(\Omega)  \to \h^1_{2\pi}$ is an
$\sone$-equivariant, compact, gradient operator given by the formula  
$$\ds \langle \nabla
\zeta(q),q_1\rangle_{\h^1_{2\pi}}= \int_0^{2\pi} \left( q(t)+ \nabla U(q(t)),q_1(t) \right) \; dt.$$
 In  other words the gradient $\nabla \Phi :
\h^1_{2\pi}(\Omega)   \rightarrow \h^1_{2\pi}$ is an $\sone$-equivariant $C^1$-operator in the form of a compact perturbation of the
identity. It is known that solutions of system \eqref{sys} are in one to one correspondence with $\sone$-orbits of solutions of 
$\nabla\Phi(q)=0.$ 
 From now on we assume that  $q_0 \in (\nabla U)^{-1}(0).$ 

Consider the linearization of the system \eqref{sys} at $q_0$ of the form 
% \beq% \label{sysl}  
$$
\left\{
\ba{rcl}  \ddot q(t) & = & - \nabla^2 U(q_0)(q-q_0) \\
q(0) & = & q(2\pi) \\
\dot q(0) & = & \dot q(2\pi)
\ea. 
\right.
%\eeq 
$$

The corresponding functional  $\Psi : \h^1_{2\pi} \rightarrow \bR$  is defined as follows

\begin{eqnarray}
\label{fun} \Psi(q) & = & \frac{1}{2} \int_0^{2\pi}  \| \dot q(t) \|^2 - (\nabla^2 U(q_0)q(t),q(t)) + 2 (\nabla^2 U(q_0)q_0,q(t))\; dt   \nonumber \\
   & = & \frac{1}{2} \|q\|^2_{\h^1_{2\pi}} - \frac{1}{2} \int_0^{2\pi}  ((\nabla^2 U(q_0)+ Id)q(t),q(t)) - 2 (\nabla^2 U(q_0)q_0,q(t))\; dt   \\ 
    & = & \frac{1}{2} \|q\|^2_{\h^1_{2\pi}} + \langle \nabla^2 U(q_0)q_0,q \rangle_{\h^1_{2\pi}} - \frac{1}{2} \int_0^{2\pi}  ((\nabla^2 U(q_0)+ Id)q(t),q(t)) \; dt  \nonumber  \\
 & = & \frac{1}{2} \|q\|^2_{\h^1_{2\pi}} + \langle \nabla^2 U(q_0)q_0,q \rangle_{\h^1_{2\pi}} - \frac{1}{2} \langle Lq,q\rangle_{\h^1_{2\pi}}. \nonumber 
\end{eqnarray}
where   $L : \h^1_{2\pi} \rightarrow \h^1_{2\pi}$ is a linear, self-adjoint,
$\sone$-equivariant and compact operator.
It is clear that $\nabla \Psi(q)=q-Lq +\nabla^2 U(q_0)q_0.$

Given $q \in \h^1_{2\pi}$ with Fourier series
$\ds q(t)= a_0   + \sum_{k =1}^{\infty} a_k \cdot   \cos   k  t   +
 b_k \cdot   \sin   k   t,$  we know that
\beq\label{fourier}
\nabla \Psi(q)=  -\nabla^2 U(q_0)( a_0-q_0)   + \sum_{k=1}^{\infty}
 (\Lambda(k) \cdot a_k)   \cdot \cos  k t +  (\Lambda(k) \cdot b_k)  \cdot \sin k t,\eeq
 where $\ds \Lambda(k)=\left(\frac{k^2}{ k^2  + 1} Id - \frac{1}{k^2 + 1} \nabla^2 U(q_0)\right)$ (see lemma 5.1.1 of \cite{[FUGORY]} for details).

Let $G$ be a compact Lie group. Denote  by $\subg$ the set of all closed subgroups of  $G.$ Two subgroups $H, H' \in \subg$ are said to be conjugate in $G$  if there is $g \in G$ such that $H=gH'g^{-1}.$ The conjugacy is an equivalence relation on $\subg.$ The class of $H \in \subg$ will be denoted by $(H)_G$ and the set of conjugacy classes will be denoted by $\subgc.$ Denote by  $\rho : G \to \on$  a continuous homomorphism. The space $\bR^n$
with the $G$-action defined by $G \times \bR^n \ni (g,x) \to \rho(g)x \in \bR^n$ is said to be a real,  orthogonal representation of $G$ which we write $\bV=(\bR^n,\rho).$  To simplify notations we write $gx$ instead of $\rho(g)x.$
 
If $x \in \bR^n$ then a group $G_x=\{g \in G : g x = x\} \in \subg$ is called the isotropy group of $x$ and  $G(x)=\{gx : g \in G\}$ is the orbit through  $x.$ Note   the orbit $G(x)$ is a smooth $G$-manifold $G$-diffeomorphic to $G \slash G_x.$  An open subset $\Omega \subset \bR^n$ is called $G$-invariant if  $G(x) \subset \Omega$ for every $x \in \Omega.$ 
  
%Fix $H \in \subg.$ Set $\Omega_H=\{x \in \Omega : G_x=H\}$ and $\Omega^H=\{x \in \Omega  : H \subset G_x\}.$ It is known that the set $\Omega_H$ is  open in $\Omega^H,$ see \cite{[DIECK],[KBO]}.    Therefore $\Omega_{\{e\}}$ is an open $G$-invariant subset of $\Omega.$

Two orthogonal representations of $G,$ say $\bV=(\bR^n,\rho), \bV'=(\bR^n,\rho'),$ are equivalent (briefly $\bV \approx_G \bV'$) if there is an equivariant linear isomorphism $L : \bV \to \bV'$ i.e. the isomorphism $L$ satisfies $L(gx)=gL(x)$ for any $g \in G, x \in \bR^n.$ 
Put $D(\bV)=\{x \in \bV : \| x \| \leq 1\}, S(\bV)=\partial D(\bV)$ and $S^{\bV}=D(\bV) \slash S(\bV).$ Since the representation $\bV$ is orthogonal, these sets are $G$ invariant. 

Denote by $\bR[1,m], m \in \bN,$   a two-dimensional representation of the group $\sone$ with an action of   $\sone$
given by $(\Phi(e^{\phi}),(x,y)) \longrightarrow (\Phi(e^{\phi}))^m(x,y)^T,$ where  $\Phi(e^{i\phi})=
\left[\begin{array}{lr}
\cos \phi & -\sin \phi\\
\sin \phi&\cos \phi
\end{array}\right].$
For $k, m \in \mathbb{N}$ we denote by $\mathbb{R}[k,m]$ the direct sum of $k$ copies of $\bR[1,m]$, we also
denote by $\mathbb{R}[k,0]$  the  $k$-dimensional trivial representation of $\sone.$ The following classical result gives a complete classification (up to an equivalence) of finite-di\-men\-sio\-nal $\sone$-representations, in \cite{[ADM]} you can find a proof of it.

\begin{Theorem}
\label{tk} If  $\bV$ is an $\sone$-representation    then there exist finite se\-quen\-ces $\{k_i\},\, \{m_i\}$
satisfying
\beq\label{rep}
 m_i\in \{0\}\cup \mathbb{N},\quad  k_i\in \mathbb{N},\quad  1\le i \le r,\,   m_1  < m_2 < \dots  <
m_r
\eeq
such that $\bV$ is equivalent to $\ds \bigoplus^r_{i=1} \bR[k_i ,m_i]$ i.e.   $\bV \approx_{\sone} \ds \bigoplus^r_{i=1} \bR[k_i,m_i]$. Moreover, the
equivalence class of $\bV$  is uniquely determined
by sequences $\{k_i\}, \{m_i\} $ satisfying \ref{rep}.
\end{Theorem}

Assume that $H \in \subg.$ Let $\bY$ be a $H$-space. The product $G \times \bY$ carries $H$-action $(h,(g,y)) \to (gh^{-1},hy).$ The orbit space of $H$-action is denoted by $G \times_H \bY$ and called the twisted product over $H.$ $G \times_H \bY$ is a $G$-space with $G$-action defined by $(g',[g,x]) \to [g'g, x].$

Let $\bY$ be a pointed $H$-space with a base point $\ast.$ Denote by  $G^+$  the group $G$ with  disjoint $G$-fixed base point $\ast$ added. Define the smash product of $G^+$ and $\bY$ by $G^+ \wedge \bY=G^+ \times \bY \slash G^+ \vee \bY =G \times \bY \slash G \times \{\ast\}.$ The group $H$ acts on the pointed space $G^+ \wedge \bY$ by $(h,[g,y]) \to [gh^{-1},hy].$ The orbit space is denoted by $G^+ \wedge_H \bY$ and called the smash over $H,$ see \cite{[DIECK]}. A~formula  $(g',[g,y]) \to [g'g,y]$ induces $G$-action so that $G^+ \wedge_H \bY$ becomes a pointed $G$-space.

Below we introduce the notion of an admissible pair of compact Lie groups. Such pairs will play crucial role in computations of the equivariant Conley index of a non-degenerate orbit of critical points of invariant potentials. 

\bdf \label{admp}
Fix $H \in \subg.$ A pair $(G,H)$ is said to be \textit{admissible} if for any $K_1,K_2 \in \subh$  the following condition is satisfied: 
$\text{ if } (K_1)_H \neq (K_2)_H \text{ then } (K_1)_G \neq (K_2)_G.$
\edf

\br \label{subg}
Of course if $(K_1)_H = (K_2)_H$ then $(K_1)_G = (K_2)_G.$ Therefore a pair $(G,H)$ is admissible if for any $K_1, K_2 \in \subh$ the condition $(K_1)_H = (K_2)_H$ is equivalent to  $(K_1)_G = (K_2)_G.$
Let $\widehat{G}$ be a compact Lie group such that $G \in \sub(\widehat{G})$  and $H \in \subg.$ If the pair $(G,H)$ is not admissible then the pair $(\widehat{G},H)$ is not admissible.
\er

\bex
Let $G=SO(4)$ and $H= \sotwo \times \sotwo.$ We claim that the pair $(G,H)$ is not admissible. Indeed, define $g=g^{-1} =\left[\ba{cc} \Theta & Id_2 \\Id_2 & \Theta  \ea \right] \in SO(4)$ where $\Theta \in \sotwo$, we observe that for $K_1=\{e\} \times \sotwo, K_2=\sotwo \times \{e\} \in \subh$ we obtain 
$$(K_1)_H = K_1 \neq K_2 =(K_2)_H \text{ and } g K_1 g^{-1}=K_2 \text{ i.e. } (K_1)_G = (K_2)_G.$$ By remark \ref{subg} we obtain that the pair $(SO(n),H), n \geq 4,$ is not admissible.
\eex

\bex
Let $\bT \in \sub (\sot)$ be the group of symmetries of tetrahedron i.e. $\bT = \{\operatorname{id}, (123),$ $ (132),$ $  (124), (142), (134), (143), (234), (243), (12)(34), (13)(24), (14)(23)\}$ is the group of even permutations of the set $\{1,2,3,4\}$ and define the group  $H = \{\operatorname{id},\; (12)(34),\; (13)(24),\;$ $ (14)(23)\} \in \subt.$ Note that $H$ is commutative. We claim that the pair $(\bT,H)$ is not admissible. Indeed, fix $g=(123) \in \bT$ and note that for  $K_1=\{\operatorname{id},(12)(34)\}, K_2 = \{\operatorname{id}, (13)(24)\}$ $ \in \subh$   we obtain 
$(K_1)_H=K_1 \neq K_2=(K_2)_H \text{ and } gK_2g^{-1}=K_1 \text{ i.e. } (K_1)_{\bT}=(K_2)_{\bT};$
which proves that the pair $(\bT,H)$ is not admissible. Moreover, by remark \ref{subg} we obtain that the pair $(SO(n),H), n \geq 3,$ is not admissible.
 \eex  

\br 
If $G$ is commutative then for any $H \in \subg$ the pair $(G,H)$ is admissible. 
\er

\bl \label{gh}
If $\Gamma$ is a compact Lie group, $G = \Gamma \times \sone$ and $H=\{e\} \times \sone$ then the pair $(G,H)$ is admissible.
\el
\begin{proof} 
For all $K_1, K_2 \in \subh$ we have  
$(K_1)_H \neq (K_2)_H \text{ iff } K_1 \neq K_2 \text{ iff } \card K_1 \neq \card K_2.$ Finally, if $ \card K_1 \neq \card K_2 \text{ then } (K_1)_G \neq (K_2)_G,$ which completes the proof.
\end{proof} 

\bl 
If $H \in \subotwo,$ then the pair $(\sot,H)$ is admissible.
\el
\begin{proof}
 Suppose, contrary to our claim, that there  is $H \in \subotwo$ and $K_1, K_2 \in H$ such that $(K_1)_H \neq (K_2)_H$ and $(K_1)_{\sot} = (K_2)_{\sot}$
By theorem 6.1 of \cite{[GSS]} every planar subgroup of $\sot$ is conjugate in $\sot$ to one of $\{e\}, \bZ_n (n \geq 2), D_n (n \geq 2), $ $ \sotwo, \otwo \in \subotwo.$ Since the cardinalities of adjoint groups are equal and  $(\sotwo)_{\sot} \neq (\otwo)_{\sot},$ equality $(K_1)_{\sot} = (K_2)_{\sot}$ implies $\card K_1 = \card K_2 < \infty.$ Taking into account that $K_1 \neq K_2$ we obtain that there is $n \geq 2$ such that $K_1=\bZ_{2n}$ and $K_2=D_n.$
Finally we obtain $(\bZ_{2n})_{\sot} = (D_n)_{\sot},$ a contradiction. 
\end{proof}

We denote by $\cF_{\ast}(G)$ the set of 
finite pointed $G$-CW-complexes and by $\cF_{\ast}[G]$ the set of
$G$-homotopy types of elements of $\cF_{\ast}(G)$.   Note that $S^{\bV} \in \cF_{\ast}(G).$ By $[\bX]_G \in \cF_{\ast}[G]$ we denote the $G$-homotopy class of $\bX \in \cF_{\ast}(G).$  Let $\mathbf{F}$ be the free abelian group generated by the elements of $\cF_{\ast}[G]$ and let $\mathbf{N}$ be the  subgroup of $\mathbf{F}$ generated by
all elements   $[\bA] - [\bX] +[\bX \slash \bA]$ for pointed $G$-CW-subcomplexes $\bA$ of a pointed $G$-CW-complex $\bX$.

\bdf \label{tdr} Let be $U(G) = \mathbf{F} \slash \mathbf{N}$ and let $\chi_G(\bX) \in U(G)$ be the class of $[\bX]$ in
$U(G)$. The element $\chi_G(\bX)$ is said to  be a $G$-equivariant  Euler  characteristic of a pointed
$G$-CW-complex $\bX$. \edf

\nt For  $\bX, \bY \in \cF_{\ast}(G)$ let $[\bX \vee \bY] \in \cF_{\ast}[G]$ denote a $G$-homotopy type of
the wedge $\bX \vee \bY \in \cF_{\ast}(G)$. Since $ [\bX] -
[\bX \vee \bY] + [\bY] = [\bX] - [\bX \vee \bY] + [(\bX \vee \bY) \slash \bX]  \in \bN$, the sum is well-defined
\begin{equation}
\label{add} \chi_G(\bX) + \chi_G(\bY) = \chi_G(\bX \vee \bY).
\end{equation}

\noindent For  $\bX, \bY \in \cF_{\ast}(G)$ let $\bX \wedge \bY = \bX \times \bY \slash \bX \vee \bY$. The
assignment $(\bX,\bY) \rightarrow \bX \wedge \bY$ induces a product $U(G) \times U(G) \rightarrow U(G)$ given by
\begin{equation}
\label{mul} \chi_G(\bX) \star \chi_G(\bY) = \chi_G(\bX \wedge \bY).
\end{equation}

\bl[\cite{[DIECK]}] \label{ms} $\left(U(G), +, \star\right)$ with an additive and multiplicative structures given
by (\ref{add}), (\ref{mul}), respectively, is a  commutative ring with unit $\bI=\chi_G(G \slash G^+).$ \el

We call $\left(U(G), +, \star\right)$   the Euler ring of a compact Lie group G.

\bl[\cite{[DIECK]}] \label{pc} The Euler  ring $\left(U(G), +, \star\right)$ is the free abelian group with
basis $\displaystyle \chi_G\left(G\slash H^+\right)$, where $(H) \in \subgc$.  Moreover, if $\bX \in \mathcal{
F}_{\ast}(G)$  then 
\beq \label{euler} \ds  \chi_G(\bX) = \sum_{(K)_G \in
\subgc} n^G_{(K)_G}(\bX) \cdot \chi_G\left(G\slash K^+\right), 
\eeq  where $\ds n^G_{(K)_G}(\bX)=\sum_{i=0}^{\infty} (-1)^{i} n(\bX,(K)_G,i)$ and $n(\bX,(K)_G,i)$ is the number of $i$-cells of type $(K)_G$ of $\bX.$ \el

Here and subsequently $\chig : \cF_{\ast}[G] \to U(G)$ stands for the equivariant Euler characteristic for finite pointed $G$-CW-complexes, see properties (IV.1.5) of \cite{[DIECK]}. 

\br \label{invert}
The Euler ring $U(\sone)$   is generated by elements $\bI=\chi_{\sone}({\sone \slash \sone}^+), \chi_{\sone}({\sone \slash \bZ_k}^+) \in U(\sone), k \in \bN.$ Since 
\beq\label{mulsone} 
\chi_{\sone}({\sone \slash \bZ_k}^+) \star  \chi_{\sone}({\sone \slash \bZ_{k'}}^+) =\Theta \in U(\sone), \eeq 
for $k, k' \in \bN,$ it is easy to see that if the  representation $\bW$ of $\sone$ is equivalent to the representation $\bR[k_0,0] \oplus \bR[k_1,m_1] \oplus \ldots \oplus \bR[k_r,m_r]$ ($\bW \approx_{\sone} \bR[k_0,0] \oplus \bR[k_1,m_1] \oplus \ldots \oplus \bR[k_r,m_r]$), then 
\beq \label{ches}
\chi_{\sone}\left(S^{\bW}\right)=\chi_{\sone}\left(S^{ \bR[k_0,0] \oplus \bR[k_1,m_1] \oplus \ldots \bR[k_r,m_r]}\right)=(-1)^{k_0}\left(\bI - \sum_{i=1}^r k_i  \chi_{\sone}(\sone \slash {\bZ_{m_i}}^+)\right).
\eeq
We claim that $\chi_{\sone}\left(S^{\bW}\right)$ is invertible in the Euler ring $U(\sone).$
Indeed by formula \eqref{mulsone} and \eqref{ches}  obtain 
$$\chi_{\sone}\left(S^{\bW}\right) \star \left((-1)^{k_0}\left(\bI + \sum_{i=1}^r k_i  \chi_{\sone}(\sone \slash {\bZ_{m_i}}^+)\right)\right)=\bI - \left(\sum_{i=1}^r k_i  \chi_{\sone}(\sone \slash {\bZ_{m_i}}^+)\right)^2 = \bI,$$ which completes the proof.
\er

The principal significance of the following theorems is that they allow us to express the $G$-equi\-va\-riant Euler characteristic of a $G$-CW-complex $G^+ \wedge_H X$ in terms of a $H$-equivariant Euler characteristic of the $H$-CW-complex $X.$ Later on, using these relations,  we will distinguish equivariant Conley indices of non-degenerate orbits of critical points of invariant potentials.

\bt \label{hsm} Let be fixed $H \in \subg$   and $X$ a $H$-CW-complex. \\  
If  $\ds \chih([X]_H)=\sum_{(K)_H \in \subhc} n^H_{(K)_H}(X) \cdot \chi_H(H \slash K^+) \in U(H),$ then
\begin{enumerate}
\item  $G^+ \wedge_H X$ is a $G$-CW complex,
\item if $\ds K \in \subg$ then  $\ds n^G_{(K)_G}(G^+ \wedge_H X)=\sum_{(\cK)_H \in \subhc, (\cK)_G=(K)_G} n^H_{(\cK)_H}(X) \in \bZ$ \\
and  $\ds \chig\left([G^+ \wedge_H X]_G\right)= \sum_{(K)_G \in \subgc} n^G_{(K)_G}(G^+ \wedge_HX) \cdot \chi_G(G \slash K^+) \in U(G).$ \end{enumerate}
\et 
\begin{proof}
First of all note that \beq \label{smpro} G^+ \wedge_H X =(G \times_H X) \slash (G \times_H \{\ast\}). \eeq \\ (1) By   proposition (II.1.13) of \cite{[DIECK]} $(G \times_H X, G \times_H \{\ast\})$ is a relative $G$-CW-complex. Therefore by formula \eqref{smpro} and  exercise (II.1.17.1) of \cite{[DIECK]} we obtain $[G^+ \wedge_H X]_G \in \cF_{\ast}[G].$ \\ (2)  Let $\{(k_1,(K_1)_H),\ldots,(k_s,(K_s)_H)\}$ be a CW-decomposition of the $H$-CW-complex $X$ i.e. $X$ consists of $s$ equivariant cells and the $j$-th cell is of dimension $k_j \in \bN \cup \{0\}$ and  orbit type $(K_j)_H.$ By proposition (II.1.13) of \cite{[DIECK]} and formula \eqref{smpro} we obtain that $\{(k_1,(K_1)_G),\ldots,(k_s,(K_s)_G)\}$ is a CW-decomposition of the $G$-CW-complex $G^+ \wedge_H X$.   It can happen that  $(K_i)_H \not = (K_j)_H $ and $ (K_i)_G =(K_j)_G.$ Thus taking into account the $G$-CW-decomposition of the $G$-CW-complex $G^+ \wedge_H X$ we obtain 
$$\ds n^G_{(K)_G}(G^+ \wedge_H X)=\sum_{(K_i)_H 
 \in \subhc, (K_i)_G=(K)_G} n^H_{(K_i)_H},$$  which completes the proof.
\end{proof}

We can significantly simplify theorem \ref{hsm} assuming that the pair $(G,H)$ is admissible.
In the theorem below we consider this case. 

\bt  \label{hsmadm} We fix $H \in \subg$ in such that it satisfies that the pair $(G,H)$ is admissible and is also a $H$-CW-complex  $X.$   If  $\ds \chih([X]_H)=\sum_{(K)_H \in \subhc} n^H_{(K)_H}(X) \cdot \chi_H(H \slash K^+) \in U(H),$ then
\begin{enumerate}
\item if  $K_1, K_2 \in \subh$ and $(K_1)_H \neq (K_2)_H \in \subhc$ then $\chi_G(G \slash K_1^+) \neq \chi_G(G \slash K_2^+) \in U(G),$
\item if $K \in \subh$ then  $n^G_{(K)_G}(G^+ \wedge_H X)=n^H_{(K)_H}(X) \in \bZ,$ and \\
$\ds \chig\left([G^+ \wedge_H X]_G\right)= \sum_{(K)_H \in \subhc} n^H_{(K)_H}(X) \cdot \chi_G(G \slash K^+) \in U(G).$ \end{enumerate}
\et  
\begin{proof}   (1) Since the pair $(G,H)$ is admissible and $(K_1)_H \neq (K_2)_H,$    $(K_1)_G \neq (K_2)_G$  Thus the $G$-spaces $G \slash K_1$ and $G \slash K_2$ are not $G$-equivalent and that is why $\chi_G(G \slash K_1^+) \neq \chi_G(G \slash K_2^+) \in U(G),$  (see \cite{[DIECK]} for more details). \\ (2) Let $\{(k_1,(K_1)_H),\ldots,(k_s,(K_s)_H)\}$ be a CW-decomposition of the $H$-CW-complex $X.$  By proposition (II.1.13) of \cite{[DIECK]} and formula \eqref{smpro} we obtain that $\{(k_1,(K_1)_G),\ldots,(k_s,(K_s)_G)\}$ is a CW-decomposition of the $G$-CW-complex $G^+ \wedge_H X$.   Now from Theorem \ref{hsm} we get  
$$\ds n^G_{(K)_G}(G^+ \wedge_H X)=\sum_{(\cK)_H \in \subhc, (\cK)_G=(K)_G} n^H_{(\cK)_H}(X) \in \bZ.$$ 
Since the pair $(G,H)$ is admissible, $(K_i)_H=(K_j)_H $ iff $(K_i)_G=(K_j)_G.$
 Consequently,   we obtain    $n^G_{(K)_G}(G^+ \wedge_H X)=n^H_{(K)_H}(X),$  which completes the proof.
\end{proof} 
 
As a direct consequence of the above theorem we obtain the following corollary.
\bco \label{hsmadmco}
Fix $H \in \subg$ such that the pair $(G,H)$ is admissible.  If $[X]_H, [Y]_H  \in \cF_{\ast}[H]$ and $\ds \chih([X]_H) \neq  \chih([Y]_H) \in U(H),$ then 
$\chig\left([G^+ \wedge_H X]_G\right) \neq \chig\left([G^+ \wedge_H Y]_G\right) \in U(G).$  In other words the map  $U(H) \ni \chih([X]_H) \to \chig\left([G^+ \wedge_H X]_G\right) \in  U(G)$ is injective.
\eco 

 %%%%%%%%%%%%%%%%%%%%%%%%%%%%%%%%%%%%%%%%%%%%%%%%%%%%%%%
\numsec
\section{The equivariant Conley index}
\label{conley}

In this section we  express the equivariant Conley index of a non-degenerate critical orbit of an invariant potential  in terms of the equivariant Conley index of an isolated critical point of this   potential restricted to the space orthogonal to this orbit. This relation allows us to distinguish Conley indices of two non-degenerate critical orbits. Let be $G$ a compact Lie group. We
fix $k \in \bN \cup \{\infty\}$ and an open and $G$-invariant subset $\Omega \subset \bR^n.$

\bdf \label{potential} A map $\vp : \Omega \to \bR$ of class $C^k$ is called \textit{ $G$-invariant $C^k$-potential}, if $\vp(gx)=\vp(x)$ for every $g \in G$ and $x \in \Omega.$ The set of $G$-invariant $C^k$-potentials   will be denoted by $C^k_{G}(\Omega,\bR).$  \edf

Fix $\vp \in C^2_G(\Omega,\bR)$ and denote by $\nabla \vp, \nabla^2 \vp$ the gradient and the Hessian of $\vp,$ respectively.
For $x_0 \in \Omega \cap  (\nabla \vp)^{-1}(0)$ denote by  $\morse(\nabla^2 \vp(x_0))$ the Morse index of the  Hessian of $\vp$ at $x_0$  i.e. the sum of the  multiplicities of negative eigenvalues of the symmetric matrix $\nabla^2 \vp(x_0).$  It is important to observe that for any $x_0'\in G(x_0)$ the following equality holds  
$\morse(\nabla^2 \vp(x_0))=\morse(\nabla^2 \vp(x_0'))$ and $(G_{x_0})=(G_{x_0'}).$ In other words the Morse index of the Hessian $\nabla^2 \vp(x_0')$ and the conjugacy class of the isotropy group $G_{x_0'}$ do not depend on the choice of an element   $x_0' \in G(x_0).$ Therefore one can assign the Morse index $\morse(\nabla^2 \vp(x_0)) \in \bN \cup \{0\}$ and the conjugacy class $(G_{x_0})_G \in \subgc$ to the orbit $G(x_0) \subset (\nabla \vp)^{-1}(0).$

\bdf \label{family} A map $\psi : \Omega \to \bR^n$ of the class $C^{k-1}$ is said to be an \textit{$G$-equivariant $C^{k-1}$-map}, if $\psi(gx)=g \psi(x)$ for every $g \in G$ and  $x \in \Omega.$ The set of $G$-equivariant $C^{k-1}$-maps  will be denoted by $C^{k-1}_{G}(\Omega,\bR^n).$  \edf

\br \label{eqgrad} It is clear that if $\vp \in C^k_G(\Omega,\bR),$ then $\nabla \vp \in C^{k-1}_G(\Omega,\bR^n).$ Moreover,  if $x_0 \in (\nabla \vp)^{-1}(0),$ then $G(x_0)   \subset (\nabla \vp)^{-1}(0)$  i.e. the $G$-orbit of a critical point consists of critical points.  If $\nabla \vp(x_0)=0$ then  $\nabla \vp(\cdot)$ is fixed on $G(x_0).$ That is why $T_{x_0} G(x_0) \subset\ker \nabla^2 \vp(x_0) $ and consequently $\dim \ker \nabla^2 \vp(x_0) \geq \dim T_{x_0} G(x_0)=\dim G(x_0).$ 
\er
\bdf 
An orbit   $G(x_0) \subset (\nabla \vp)^{-1}(0)$ is said to be   \textit{non-degenerate}, if   $T_{x_0} G(x_0) =$ $\ker \nabla^2 \vp(x_0) $ or equivalently if $\dim \ker \nabla^2 \vp(x_0) = \dim T_{x_0}G(x_0).$\edf

Let $\vp \in C^2_G(\Omega,\bR)$ and $x_0 \in \Omega.$ Suppose that the orbit $G(x_0) \subset (\nabla \vp)^{-1}(0)$ is non-degenerate. By the equivariant Morse lemma, see \cite{[WAS]}, $G(x_0)$ is isolated in $(\nabla \vp)^{-1}(0).$ 

The rest of this section is dedicated to the study of the $G$-equivariant Conley index, see \cite{[BARTSCH],[FLOER],[GEBA],[SMWA]}, of the isolated invariant set $G(x_0)$ considered as  a $G$-orbit of stationary solutions of the equation $\dot x(t)=-\nabla \vp(x(t)).$
Note that since the orbit $G(x_0)$ is non-degenerate, $T_{x_0} \bR^n =T_{x_0} G(x_0) \oplus T_{x_0} G(x_0)^{\perp} =  \ker \nabla^2 \vp(x_0) \oplus T_{x_0} G(x_0)^{\perp}.$  

In the  following, for simplicity of notation, we write   $H$ instead of $G_{x_0}.$ 
Since $\bR^n$ is an orthogonal representation of $G,$ $T_{x_0} G(x_0)^{\perp}$ is an orthogonal representation of $H.$ Define $\phi = \vp_{\mid T_{x_0} G(x_0)^{\perp}} \in C^2_H(T_{x_0} G(x_0)^{\perp},\bR).$ Since   the orbit $G(x_0)$ is non-degenerate,  $\nabla \phi(x_0)=0$ and $\nabla^2 \phi(x_0)$ is non-degenerate. That is why, by the Morse lemma, $x_0$ is an  isolated critical point of the potential $\phi.$ 

Again by the non-degeneracy of $G(x_0)$ we have   

\beq 
\label{decs} T_{x_0} \bR^n = \ba{c}  \ker \nabla^2 \vp(x_0) \\ \oplus \\ T_{x_0}^{\perp} G(x_0)^H \\ \oplus \\ T_{x_0}^{\perp} G(x_0) \ominus T_{x_0}^{\perp} G(x_0)^H \ea. 
\eeq 

%\beq \label{decs} T_{x_0} \bR^n = \ba{c}  \ker \nabla^2 \vp(x_0) \\ \oplus \\ T_{x_0}^{\perp} G(x_0)^H \\ \oplus \\ T_{x_0}^{\perp} G(x_0) \ominus T_{x_0}^{\perp} G(x_0)^H \ea. 
%\eeq 

Moreover, we obtain the following decomposition of the Hessian $\nabla^2 \vp(x_0)$:
\beq \label{decm}\nabla^2 \vp(x_0)=\left[\ba{ccc}  \Theta & 0 & 0 \\ 0 & B(x_0) & 0 \\ 0 & 0 & C(x_0)\ea\right] : T_{x_0} \bR^n \to T_{x_0} \bR^n,\eeq see \cite{[GEBA]}. We observe that the matrices $B(x_0), C(x_0)$ are non-degenerate because the orbit $G(x_0)$  is non-degenerate.

For $\epsilon > 0$ we define $D_{\epsilon}(x_0)=\{x \in T^{\perp}_{x_0} G(x_0): \|x-x_0\| \leq \epsilon\}$ and $S_{\epsilon}(x_0)=\partial D_{\epsilon}(x_0).$  The sets   $S_{\epsilon}(x_0), D_{\epsilon}(x_0) \subset T^{\perp}_{x_0} G(x_0)$ are $H$-invariant because   $x_0 \in T^{\perp}_{x_0} G(x_0)^H.$ Since $\nabla^2 \phi(x_0)$ is self-adjoint, there is an orthogonal decomposition $T^{\perp}_{x_0} G(x_0)=T^{\perp}_{x_0} G(x_0)^+ \oplus T^{\perp}_{x_0} G(x_0)^-$ corresponding to the positive and negative part of the spectrum of $\nabla^2 \phi(x_0).$ Without loss of generality one can assume that $(D^+_{\epsilon}(x_0) \times D^-_{\epsilon}(x_0)) \cap (\nabla \phi)^{-1}(0) = \{x_0\},$ where  $D^{\pm}_{\epsilon}(x_0)=\{x \in T^{\perp}_{x_0} G(x_0)^{\pm} : \|x-x_0\| \leq \epsilon\}$ and $S^{\pm}_{\epsilon}(x_0)=\partial D^{\pm}_{\epsilon}(x_0).$

Let us compute the $H$-equivariant Conley index of the isolated invariant set $\{x_0\}$ considered as a stationary solution of the equation $\dot x(t)=-\nabla \phi(x(t)).$ Since $\nabla^2 \phi(x_0)$ is non-degenerate, without loss of generality, instead of equation $\dot x(t)=-\nabla \phi(x(t))$ we will consider its linearization at $x_0$ i.e. $\dot x(t)=-\nabla^2 \phi(x_0)x(t).$

Note that $(N,L)=(D^-_{\epsilon}(x_0) \times D^+_{\epsilon}(x_0), S^-_{\epsilon}(x_0) \times D^+_{\epsilon}(x_0))$ is a $H$-index pair of the isolated invariant set $\{x_0\}$. Consequently,  the $H$-equivariant Conley index   of   $\{x_0\}$, denoted by $\cih(\{x_0\},-\nabla \phi)$ is equal to 
$\cih(\{x_0\},-\nabla \phi)=([N \slash L]_H,[L])$ i.e. the $H$-equivariant Conley index is the $H$-homotopy type of the quotient $H$-space $N \slash L.$ Moreover, since $(D^-_{\epsilon}(x_0), S^-_{\epsilon}(x_0))$ is a strong $H$-deformation retract of $(N,L),$  we obtain the following equality 
\beq \label{hci}
\cC\cI_H(\{x_0\},-\nabla \phi)=([D^-_{\epsilon}(x_0) \slash S^-_{\epsilon}(x_0)]_H,[S^-_{\epsilon}(x_0)]) \in \cF_{\ast}[H].
\eeq

Let us compute the $G$-equivariant Conley index of the isolated invariant set $G(x_0)$ considered as a $G$-orbit of stationary solutions of the equation $\dot x(t)=-\nabla \vp(x(t)).$  Note that if $(\bX, \bA)$ is a pair of $H$-spaces and $(\bX_0,\bA_0)$ is a strong $H$-deformation retract of $(\bX, \bA)$, then  $(G \times_H \bX_0,G \times_H \bA_0)$ is a strong $G$-deformation retract of $(G \times_H \bX,G \times_H \bA),$ see \cite{[DIECK],[KBO]} for other properties of the twisted product over $H.$

First we express the $G$-index $(\cN,\cL)$ pair of the orbit $G(x_0)$ in terms of the twisted product over $H$    of the $H$-index pair $(N,L)$ of $x_0.$
In fact $(\cN,\cL)=(G \times_H N,G \times_H L)$ is a $G$-index pair of the isolated invariant set $G(x_0).$ 

Therefore the $G$-equivariant Conley index  of the isolated invariant set $G(x_0)$ is equal to $\cig(G(x_0),-\nabla \vp)=([\cN \slash \cL]_G,[\cL]) \in \cF_{\ast}[G].$
Moreover,   since $(D^-_{\epsilon}(x_0), S^-_{\epsilon}(x_0))$ is a strong $H$-deformation retract of $(N,L),$ the pair  $(G \times_H D^-_{\epsilon}(x_0) , G  \times_H S^-_{\epsilon}(x_0))$  is a strong $G$-deformation retract of $(G \times_H N,G \times_H L).$ Consequently we obtain the  following equality 
\beq 
\label{gci} \cC\cI_G(G(x_0),-\nabla \vp)= ([(G \times_H D^-_{\epsilon}(x_0))\slash (G  \times_H S^-_{\epsilon}(x_0))]_G,[G \times_HS^-_{\epsilon}(x_0)]) \in \cF_{\ast}[G].
\eeq

In the theorem below we present the relation between the equivariant Conley indices given by formulas \eqref{hci} and \eqref{gci} i.e. we express the $G$-equivariant Conley index of a  non-degenerate orbit $G(x_0) \subset (\nabla \vp )^{-1}(0)$  in terms of the $H$-equivariant Conley index of a non-degenerate  critical point $x_0 \in (\nabla \phi)^{-1}(0)$ of the restricted functional $\phi.$ 

\bt \label{cio} Let $\vp \in C^2_G(\Omega,\bR)$ and $x_0 \in \Omega.$ Suppose that the orbit $G(x_0) \subset (\nabla \vp)^{-1}(0)$ is  non-degenerate. Then
$\ci_G(G(x_0),-\nabla \vp)= G^+ \wedge_H\ci_H  (\{x_0\},-\nabla \phi) \in \cF_{\ast}[G],$ where $H= G_{x_0}.$\et
\begin{proof}
To simplify notation set $(D^-,S^-)=(D^-_{\epsilon}(x_0),S^-_{\epsilon}(x_0)).$
Note that $(D^-,S^-)$ is a relative $H$-CW-complex and that $D^- \slash S^-$ is a pointed $H$-CW-complex, see \cite{[DIECK],[MAY]}.
It is clear that  the spaces
$(G \times D^-) \slash (G \times S^-)$ and $(G \times (D^- \slash S^-)) \slash (G \times \{\ast\})=G^+ \wedge (D^- \slash S^-)$ are homeomorphic.
Consequently, taking into account the  $H$-action on $G^+ \wedge (D^- \slash S^-)$ and $G \times D^-$ we obtain that 
$([G \times_H D^- \slash G \times_H S^-]_G,[G \times_H S^-]) = ([G^+ \wedge_H (D^- \slash S^-)]_G,[\ast])=(G^+ \wedge [D^-\slash S^-]_H,\ast),$ which completes the proof.
\end{proof}

By theorem \ref{cio}  the computation of the   $G$-equivariant Conley index $\ci_G(G(x_0),-\nabla \vp) \in \cF_{\ast}[G]$ of the orbit  $G(x_0)$  can be reduced to  the computation of the  $H$-equivariant Conley index  $ \cih(\{x_0\},-\nabla \phi) \in \cF_{\star}[H]$ of the non-degenerate critical point   $x_0$ of the potential $\phi.$  

In the following corollary   we show how to distinguish $G$-equivariant Conley indices of non-degenerate orbits of critical points of a potential $\vp \in C^2_G(\Omega, \bR)$ considering the restrictions of this potential to orthogonal spaces  to these orbits. 
 
\bco \label{cio1}
Let $G(x'_0), G(x''_0) \subset  (\nabla \vp)^{-1}(0)$ be non-degenerate orbits  of the potential $\vp \in C^2_G(\Omega, \bR).$  We have
 $$\chihn (\cihn (\{x_0^{\nu}\},-\nabla \phi^{\nu})) =$$ $$=\sum_{(K^{\nu})_{H^{\nu}} \in \sub[H^{\nu}]} n^{H^{\nu}}_{(K^{\nu})_{H^{\nu}}}\left(\cihn (\{x_0^{\nu}\},-\nabla \phi^{\nu})\right) \cdot \chihn ({H^{\nu} \slash K^{\nu}}^+) \in U(H^{\nu}),$$ where $\nu= ' $ or $\nu= ''$  depending the case, $\phi^{\nu} = \vp_{\mid T^{\perp}_{x_0^{\nu}} G(x_0^{\nu})}$ and $H^{\nu}=G_{x_0^{\nu}}.$  Assume that  there is $(K_0^{\nu})_{H^{\nu}} \in \sub[H^{\nu}]$   such that
\begin{enumerate}
\item $(K'_0)_G=(K''_0)_G \in \subgc,$
\item $\ds \sum_{\sub[H'] \ni  (K)_{H'}:(K)_G =(K'_0)_G} n^{H'}_{(K)_{H'}}(\ci_{H'}(\{x_0'\},-\nabla \phi')) \neq  $ \\ $\ds \sum_{\sub[H''] \ni  (K)_{H''}:(K)_G =(K''_0)_G} n^{H''}_{(K)_{H''}}(\ci_{H''}(\{x_0''\},-\nabla \phi'')).$
\end{enumerate}
Then
$\cig(G(x'_0),-\nabla \vp) \neq  \cig(G(x''_0),-\nabla \vp) \in \cF_{\ast}[G].$
Moreover, $\chig(\cC\cI_G(G(x'_0),-\nabla \vp)) \neq  $ $\chig(\cig(G(x''_0),-\nabla \vp)) \in U(G).$
\eco
\begin{proof}  Fix $\nu \in \{',''\}$, then
by theorem \ref{cio} we obtain $\ci_G(G(x^{\nu}_0),-\nabla \vp)=G^+ \wedge_{H^{\nu}}\ci_{H^{\nu}}  (\{x_0\},-\nabla \phi^{\nu}) $ $\in \cF_{\ast}[G].$ Hence $\chi_G(\ci_G(G(x^{\nu}_0),-\nabla \vp))= \chi_G(G^+ \wedge_{H^{\nu}}\ci_{H^{\nu}}  (\{x_0\},-\nabla \phi^{\nu})) \in U(G).$ From theorem \ref{hsm} we have $$n^G_{(K^{\nu}_{0})_G}(\cig(G(x_0^{\nu}),-\nabla \vp)) = \sum_{\sub[H^{\nu}] \ni  (K)_{H^{\nu}}:(K)_G =(K^{\nu}_{0})_G} n^{H^{\nu}}_{(K)_{H^{\nu}}}(\ci_{H^{\nu}}(\{x_0^{\nu}\},-\nabla \phi^{\nu})).$$
By the above formula and assumption (2) we obtain $$n^G_{(K'_{0})_G}(\cig(G(x_0'),-\nabla \vp))  \neq n^G_{(K''_{0})_G}(\cig(G(x_0''),-\nabla \vp)).$$ Hence from assumption (1) it follows that   $\chig(\cig(G(x'_0),-\nabla \vp)) \neq  \chig(\cig(G(x''_0),-\nabla \vp)),$ which completes the proof. 
\end{proof}

In the corollary below we assume that the pair $(G,H^{\nu}),\nu \in \{',''\},$ is admissible. It allows us to control the relation between the equivariant Euler characteristics $\chihn(\cihn(\{x^{\nu}_0\},-\nabla \phi^{\nu}))  \in U(H^{\nu})$ and $\chig(\cig(G(x^{\nu}_0),-\nabla \vp))   \in U(G).$

\bco \label{hgg} Let   $G(x'_0), G(x''_0) \subset  (\nabla \vp)^{-1}(0)$ be non-degenerate orbits of critical points of the potential $\vp \in C^2_G(\Omega, \bR)$ s.t.  $G_{x'_0}=G_{x''_0}(=H).$ If the pair $(G,H)$ is admissible and
$\chi_H(\cC\cI_H(\{x'_0\},-\nabla \phi')) \neq \chi_H(\cC\cI_H(\{x''_0\},-\nabla \phi'')) \in U(H)$ then 
$$\ci_G(G(x'_0),-\nabla \vp) \neq  \cC\cI_G(G(x''_0),-\nabla \vp) \in \cF_{\ast}[G].$$
Moreover, $\chi_G(\cC\cI_G(G(x'_0),-\nabla \vp)) \neq  \chi_G(\cC\cI_G(G(x''_0),-\nabla \vp)) \in U(G).$ 
\eco
\begin{proof}  Fix $\nu \in \{',''\}$ and put in Theorem \ref{hsmadm} $[X]_H=\cih(\{x^{\nu}_0\},-\nabla \phi^{\nu}).$ Applying Theorem \ref{hsmadm} we obtain that if  
$$\chih(\cih(\{x^{\nu}_0\},-\nabla \phi^{\nu}))=\sum_{(K)_H \in \subhc } n^H_{(K)_H}(\cC\cI_H(\{x^{\nu}_0\},-\nabla \phi^{\nu})) \cdot \chih(H \slash K^+) \in U(H),$$
then 
$$ \chi_G(\cig(G(x^{\nu}_0),-\nabla \vp))=\sum_{(K)_H \in \subhc}  n^H_{(K)_H}(\cC\cI_H(\{x^{\nu}_0\},-\nabla \phi^{\nu})) \cdot \chig(G \slash K^+) \in U(G).$$ 
It follows that the assumption $\chi_H(\cC\cI_H(\{x'_0\},-\nabla \phi')) \neq \chi_H(\cC\cI_H(\{x''_0\},-\nabla \phi'')) \in U(H)$ implies  that $\chi_G(\cC\cI_G(G(x'_0),-\nabla \vp)) \neq  \chi_G(\cC\cI_G(G(x''_0),-\nabla \vp)).$ Consequently $\cig(G(x'_0),-\nabla \vp) \neq  \cig(G(x''_0),-\nabla \vp),$ which completes the proof.
 \end{proof}

The following corollary is a consequence of  Corollary \ref{hgg}. The point of the corollary  is that its assumptions are expressed in terms of the representation theory of compact Lie groups.  

\bco  \label{hggr}  
Under the assumptions of Corollary  \ref{hgg}. If moreover,
\begin{enumerate}
\item $H$ is nontrivial and connected, $\dim T^{\perp}_{x'_0} G(x'_0)^-=\dim T^{\perp}_{x''_0} G(x''_0)^-\text{ and }   T^{\perp}_{x'_0} G(x'_0)^- \not \approx_H T^{\perp}_{x''_0} G(x''_0)^-,$ or 
\item  $H$ is connected, $\dim T^{\perp}_{x'_0} G(x'_0)^- > \dim T^{\perp}_{x''_0} G(x''_0)^-  \text{ and } T^{\perp}_{x'_0} G(x'_0)^- \ominus T^{\perp}_{x''_0} G(x''_0)^- \not \approx_H \bR[2k] \: \forall k \in \bN$ i.e. the representation $T^{\perp}_{x'_0} G(x'_0)^- \ominus T^{\perp}_{x''_0} G(x''_0)^-$  of $H$  is not equivalent to an even-dimensional trivial representation of $H$
\end{enumerate}
then
$\cig(G(x'_0),-\nabla \vp) \neq  \cig(G(x''_0),-\nabla \vp) \in \cF_{\ast}[G].$
Moreover, $\chig(\cig(G(x'_0),-\nabla \vp)) \neq  \chig(\cig(G(x''_0),-\nabla \vp)) \in U(G).$
\eco 
\begin{proof} 
(1) Applying Theorem 3.2 of \cite{[RYBICKIANS]} we obtain $$\chih(\cih(\{x_0'\},-\nabla \phi')) = \chih([D'^- \slash S'^-]_H) \neq  $$ $$\chi_H([D''^- \slash S''^-]_H)=\chih(\cih(\{x_0''\},-\nabla \phi'')) \in U(H).$$ The rest of the proof is a direct consequence of Corollary \ref{hgg}.\\ (2) The proof of this case is literally the same as the proof of the previous one with Theorem 3.2  of \cite{[RYBICKIANS]} replaced by Theorem 3.1 of \cite{[RYBICKIANS]}.
\end{proof}

In the following theorem we consider special orbits $G(x_0'), G(x_0'') \subset (\nabla \vp)^{-1}(0)$  i.e. non-degenerate orbits satisfying additional assumption $\morse(C(x_0^{\nu}))=0, \nu \in \{',''\},$ see formula \eqref{decm}. It allows us to compute the CW-decompositions of a $H^{\nu}$-CW-complex $\cihn(\{x_0^{\nu}\},-\nabla \phi^{\nu}) \in \cF_{\ast}[H^{\nu}]$ and a $G$-CW-complex $\cig(G(x^{\nu}_0),-\nabla \vp) \in \cF_{\ast}[G].$

\bt \label{last} 
Let $G(x'_0), G(x''_0) \subset  (\nabla \vp)^{-1}(0)$ be non-degenerate orbits of critical points of the potential $\vp \in C^2_G(\Omega, \bR)$ such that   $\morse(C(x_0')) = \morse(C(x_0''))=0.$ For $H'=G_{x_0'}$ and $H''=G_{x_0''}$
\begin{enumerate}
\item if $(H')_G \neq (H'')_G$ then $\cig(G(x'_0),-\nabla \vp) \neq  \cig(G(x''_0),-\nabla \vp) \in \cF_{\ast}[G]$ and moreover $\chig(\cig(G(x'_0),-\nabla \vp)) \neq  \chig(\cig(G(x''_0),-\nabla \vp)) \in U(G).$
\item  if  $\morse(B(x_0')) \neq \morse(B(x_0''))$  then   $\cig(G(x'_0),-\nabla \vp) \neq  \cig(G(x''_0),-\nabla \vp) \in \cF_{\ast}[G].$
\item  if $\morse(B(x_0')) - \morse(B(x_0''))$ is odd then  $\cig(G(x'_0),-\nabla \vp) \neq  \cig(G(x''_0),-\nabla \vp) \in \cF_{\ast}[G]$ and moreover $\chi_G(\cC\cI_G(G(x'_0),-\nabla \vp)) \neq  \chi_G(\cC\cI_G(G(x''_0),-\nabla \vp)) \in U(G).$
\end{enumerate}
\et
\begin{proof} For simplicity of notations set $(D^{\nu -},S^{\nu-})= (D^-_{\epsilon}(x^{\nu}_0),S^-_{\epsilon}(x^{\nu}_0)),$ where $\nu \in \{',''\}.$  Taking into account decompositions \eqref{decs}, \eqref{decm} and formula \eqref{hci} we obtain for $\nu \in \{',''\}$ that  $\cC\cI_{H^{\nu}}(\{x_0^{\nu}\},-\nabla \phi^{\nu}) \in \cF_{\ast}[H^{\nu}]$ is a $H^{\nu}$-homotopy type of pointed $H^{\nu}$-CW-complex $X^{\nu}=([D^{\nu-} \slash S^{\nu-}]_{H^{\nu}},[S^{\nu-}]).$  The $H^{\nu}$-CW-complex $X^{\nu}$ consists of base point   and one equivariant cell of dimension $\morse(B(x_0^{\nu}))$ and orbit type $(H^{\nu})_{H^{\nu}}.$ By Theorem \ref{cio} and Proposition (II.1.13) of \cite{[DIECK]} the equivariant Conley index $\cC\cI_G(G(x^{\nu}_0),-\nabla \vp) \in \cF_{\ast}[G]$ is a $G$-homotopy type of a pointed $G$-CW-complex $Y^{\nu}=G \wedge_{H^{\nu}} X^{\nu}$ which consists of base point and one equivariant cell of dimension $\morse(B(x_0^{\nu}))$ and orbit type $(H^{\nu})_G.$  
Taking into account the above remarks and formula \eqref{euler} we obtain 
\beq \label{special} \chi_G(\cC\cI_G(G(x^{\nu}_0),-\nabla \vp))=\chi_G(Y^{\nu})=(-1)^{\morse(B(x_0^{\nu}))} \cdot \chi_G(G \slash {H^{\nu}}^+) \in U(G).
\eeq 

(1) The condition $\cC\cI_G(G(x'_0),-\nabla \vp) \neq  \cC\cI_G(G(x''_0),-\nabla \vp) \in \cF_{\ast}[G]$ is fulfilled because the only nontrivial cells of $G$-CW-complexes $Y', Y''$ are of different homotopy types $(H')_G, (H'')_G,$ respectively i.e. the $G$-CW-complexes $Y',Y''$ are not $G$-homotopically equivalent. The assumption $(H')_G \neq  (H'')_G$ implies   $\chi_G(G \slash H'^+) \neq \chi_G(G \slash H''^+) \in U(G).$ Taking into account formula \eqref{special} we complete the proof.

(2) The condition $\cC\cI_G(G(x'_0),-\nabla \vp) \neq  \cC\cI_G(G(x''_0),-\nabla \vp) \in \cF_{\ast}[G]$ is fulfilled because the only nontrivial cells of $G$-CW-complexes $Y', Y''$ are of different dimensions  $\morse(B(x_0')),  \morse(B(x_0'')),$ respectively i.e. the $G$-CW-complexes $Y',Y''$ are not $G$-homotopically equivalent.

(3) If $\morse(B(x_0')) - \morse(B(x_0''))$ is odd then,  the first part of the  assertion proved in (2), it follows that  $\morse(B(x_0')) \neq \morse(B(x_0'')).$  To simplify the argument, without loss of generality, we assume that $\morse(B(x'_0))$ is even.  Hence applying formula \eqref{special}  we obtain  $$\chi_G(\cC\cI_G(G(x'_0),-\nabla \vp)) =  \chi_G(G \slash H'^+) \neq - \chi_G(G \slash H''^+) = \chi_G(\cC\cI_G(G(x''_0),-\nabla \vp)) \in U(G),$$ which completes the proof.
\end{proof}

%%%%%%%%%%%%%%%%%%%%%%%%%%%%%%%%%%%%%%%%%%%%%%%%%%%%%%%%%%%%%%

\numsec
\section{Proof of the Symmetric Liapunov center theorem}
\label{main}

In this section,   using the equivariant Conley index defined in \cite{[IZYDOREK]},  we prove  the main result of this article,  the Symmetric Liapunov center theorem for second order differential equations with symmetric potentials stated in Theorem\ref{main-theo}.

We consider $\bR^n$ as an orthogonal representation of a compact  Lie group $\Gamma$ and denote by $\rho : \Gamma \to \on$ the representation homomorphism. Denote by $\Omega \subset \bR^n$ an open and $\Gamma$-invariant subset. Fix $U \in C^2_{\Gamma}(\Omega,\bR)$ and  $q_0 \in \Omega$ a critical point of the potential $U.$   It is clear that  $\Gamma(q_0) \subset (\nabla U)^{-1}(0)$ i.e. the orbit $\Gamma(q_0)$ consists of  critical points of $U.$   In this section we study periodic solutions of system \eqref{newsys}
in a neighborhood of the orbit $\Gamma(q_0).$ 

 We observe that the orbit $\Gamma(q_0)$ is $\Gamma$-homeomorphic to $\Gamma \slash \Gamma_{q_0}=\Gamma$, for this reason it can happen that elements of this orbit are not isolated. For example if $\Gamma=\sotwo$ acts freely on $\Omega$ then the orbit $\Gamma(q_0)$ is $\sotwo$-homeomorphic to $\Gamma \slash \Gamma_{q_0}=\Gamma \slash \{e\} =\sotwo \approx  \sone.$ Note that if the group $\Gamma$ is trivial then we obtain the classical Liapunov center theorem, see \cite{[BERGER1],[BERGER2],[MAWI],[SZULKIN]} and references therein. 

Before we begin with the proof of  Theorem \ref{main-theo} we will prove one technical lemma.
This lemma will be the key ingredient in the proof of our main result.
Note that the study of periodic solutions of system  \eqref{newsys} of any period is equivalent to the study of $2\pi$-periodic solutions of the following family 
\beq \label{jesysl1}   
\left\{
\ba{rcl}  \ddot q(t) & = & -\lambda^2 \nabla U(q(t)) \\
q(0) & = & q(2\pi) \\
\dot q(0) & = & \dot q(2\pi)
\ea.
\right.
\eeq 
The $2 \pi \lambda$-periodic solution of  system \eqref{newsys} corresponds to $2\pi$-periodic solutions of \eqref{jesysl1}. Since   $\Gamma(q_0) \subset (\nabla U)^{-1}(0),$ for every $\lambda > 0$ the orbit $\Gamma(q_0)$ consists stationary solutions of system \eqref{jesysl1}.
Periodic solutions of system \eqref{jesysl1} can be considered as critical orbits of  $G=(\Gamma \times \sone)$-invariant potential of the class $C^2$ defined on $\h^1_{2\pi}(\Omega) \times (0,+\infty).$
 It is easy to show
that $\left(\h^1_{2\pi},\langle\cdot,\cdot\rangle_{\h^1_{2\pi}}\right)$ is an orthogonal representation of the group $G$
with a $G$-action defined by 
$$G \times \h^1_{2\pi} \ni ((\gamma,e^{i \theta}),q(t)) \to \gamma q(t+\theta) \text{ mod } 2 \pi.$$ It is clear tat $\h^1_{2\pi}(\Omega)$ is open and $G$-invariant.
Moreover, $\bH^1_{2\pi} = \overline{\h_0 \oplus \bigoplus_{k=1}^{\infty} \h_k} $ 
where $\h_k, k=0,1,\ldots$ are orthogonal representations of $G,$ see formula \eqref{space} and the text below it.
Define a $G$-invariant  functional 
$\Phi : \h^1_{2\pi}(\Omega) \times (0,+\infty) \to \bR$ of the class $C^2$ by  
$
 \ds \Phi(q,\lambda) = \int_0^{2\pi} \frac{1}{2} \| \dot q(t) \|^2 - \lambda^2 U(q(t)) \; dt.
$

It is well known that the solutions of  equation 
\beq \label{gradf} \nabla_q \Phi(q,\lambda)=0,\eeq are in one to one correspondence with the solutions of  system  \eqref{jesysl1}.  Since $q_0 \in \h^1_{2\pi}(\Omega)$ is a constant function, $G(q_0)=\Gamma(q_0) \subset \bH_0 = \bR^n \subset \h^1_{2\pi}$ solve equation \eqref{gradf} for any $\lambda > 0.$ 
The set of solutions $\cT=G(q_0) \times (0,+\infty) \subset \h^1_{2\pi}(\Omega) \times (0,+\infty)$ of equation \eqref{gradf} we treat as a  family of trivial solutions of equation \eqref{gradf}. 

To   proof Theorem \ref{main-theo} we will study solutions of equation \eqref{gradf}. More precisely, we will apply theorems of equivariant bifurcation theory to prove the existence of local bifurcation of solutions of equation \eqref{gradf} from the family $\cT.$ As a topological tool we will use the equivariant Conley index, see \cite{[BARTSCH],[FLOER],[GEBA],[SMWA]} and its infinite-dimensional generalization \cite{[IZYDOREK]}.

We claim that bifurcation of solutions of equation \eqref{gradf} from the trivial family $\cT$ can occur only at degenerate levels. Below we characterize 
these levels.

Indeed, since for every $\lambda > 0$ the gradient $\nabla_q \Phi(\cdot, \lambda)$ is constant on the orbit $G(q_0) \subset \h^1_{2\pi}(\Omega),$  $ \dim \ker \nabla_q^2 \Phi(q_0,\lambda) \geq \dim G(q_0).$ Applying the equivariant implicit function theorem  one can prove that bifurcation can occur   only at degenerate orbits 
$G(q_0) \times \{\lambda_0\} \subset \cT$ i.e. orbits satisfying the following condition 
\beq \label{ncl}  \dim \ker \nabla_q^2 \Phi(q_0,\lambda_0) > \dim G(q_0).\eeq 

To find parameters $\lambda_0$ satisfying condition \eqref{ncl} we   study $\dim \ker \nabla^2_q \Phi(q_0,\lambda)$ along the family of trivial orbits $\cT.$
It is clear that the study of $\ker \nabla^2_q \Phi(q_0,\lambda)$ is equivalent to  study solutions of system  
%\beq \label{sysl}  
$$
\left\{
\ba{rcl}  \ddot q(t) & = & - \lambda^2 \nabla^2 U(q_0)(q-q_0) \\
q(0) & = & q(2\pi) \\
\dot q(0) & = & \dot q(2\pi)
\ea.
\right.
%\eeq 
$$

Without loss of generality we can assume that $J(\nabla^2 U(q_0)) = \nabla^2 U(q_0)$  i.e. $\nabla^2 U(q_0)$ is in the Jordan normal form. Since the Hessian $\nabla^2 U(q_0)$ is symmetric, it is diagonal. This assumption will simplify arguments in the rest of the proof without loss of generality.

Using equation \eqref{fourier} we obtain that condition \eqref{ncl}  is fulfilled iff $\lambda \in \Lambda = \{k\slash \beta_j: k \in \bN \text{ and } j=1,\ldots,m\}.$ In other words bifurcations of $G$-orbits of solutions of equation \eqref{gradf}   from the trivial family $\cT$ can occur 
only at orbits $G(q_0) \times \Lambda \subset \cT.$ 

Fix $\beta_{j_0}$ satisfying the assumptions of Theorem \ref{main-theo}, choose  $\epsilon > 0$ sufficiently small and define $\ds \lambda_{\pm}=\frac{1\pm \epsilon}{\beta_{j_0}}.$ Without loss of generality one can assume that $[\lambda_-,\lambda_+] \cap \Lambda = \{1/\beta_{j_0}\}.$ 

Since  $G(q_0) \subset \h^1_{2\pi}$  is a non-degenerate critical orbit of the $G$-invariant functional $\Phi(\cdot,\lambda_{\pm}) : \h^1_{2\pi} \to \bR,$  it is isolated  in $\nabla \Phi(\cdot,\lambda_{\pm})^{-1}(0).$ Therefore $G(q_0)$ is an isolated invariant set in the sense of the $G$-equivariant Conley index theory defined in \cite{[IZYDOREK]} i.e. $\cig(G(q_0),-\nabla \Phi(\cdot,\lambda_{\pm}))$ is defined. 

Let $\h \subset \h^1_{2\pi}$ be a linear subspace normal    to $G(q_0)$ at $q_0$ i.e. $\h=T_{q_0}^{\perp} G(q_0) \subset \h^1_{2\pi}.$ Since the isotropy group $\Gamma_{q_0}$ is trivial, $G_{q_0}=\{e\} \times \sone$ and $\h \:$ is an orthogonal representation of $\sone.$ Define an $\sone$-invariant functional of the class  $C^2$ by $\Psi_{\lambda_{\pm}}=\Phi(\cdot,\lambda_{\pm})_{\mid \bH} : \bH \to \bR.$ Since $G(q_0) \subset \bH^1_{2\pi}(\Omega)$  is a non-degenerate critical orbit of the $G$-invariant functional $\Phi(\cdot,\lambda_{\pm}),$  $q_0 \in \bH$ is a non-degenerate critical point of $\sone$-invariant potential $\Psi_{\lambda}.$  Hence $q_0$ is an isolated invariant set in the sense of the $\sone$-equivariant Conley index theory defined in \cite{[IZYDOREK]} i.e. $\ci_{\sone}(\{q_0\},-\nabla \Psi_{\lambda_{\pm}})$ is defined. 

Note that $\bH_0=\bR^n= T^{\perp}_{q_0}\Gamma(q_0) \oplus  T_{q_0}\Gamma(q_0)$ and  $\ds \bH^n= T_{q_0}^{\perp} \Gamma(q_0) \oplus \bigoplus_{k=1}^{n} \bH_k \subset \bH$. We define $\Psi^n_{\lambda_{\pm}}={\Psi_{\lambda_{\pm}}}_{\mid \bH^n} : \bH^n \to \bR$ and note that $q_0 \in \bH^n$ is a non-degenerate critical point of a $\sone$-invariant potential $\Psi^n_{\lambda_{\pm}}.$

\bl  \label{cis} There exists $n_0 \in \bN$ such that  for any $n \geq n_0$ $$\chi_{\sone}(\ci_{\sone}(\{q_0\},-\nabla \Psi^n_{\lambda_-})) = \chi_{\sone}(\ci_{\sone}(\{q_0\},-\nabla \Psi^{n_0}_{\lambda_-})) \neq $$ $$\chi_{\sone}(\ci_{\sone}(\{q_0\},-\nabla \Psi^{n_0}_{\lambda_+}))= \chi_{\sone}(\ci_{\sone}(\{q_0\},-\nabla \Psi^n_{\lambda_+})).$$ Moreover,  $\ci_{\sone}(\{q_0\},-\nabla \Psi_{\lambda_-}) \neq \ci_{\sone}(\{q_0\},-\nabla \Psi_{\lambda_+}).$ \el
\begin{proof} 
Fix $n \in \bN.$
Since   $G(q_0)$ is non-degenerate orbit of the $\Phi(\cdot,\lambda_{\pm})$,  $\nabla \Psi^n_{\lambda_{\pm}}(q_0)=0$ and $\nabla^2 \Psi^n_{\lambda_{\pm}}(q_0)$ is an isomorphism. Consequently we obtain $\nabla \Psi^n_{\lambda_{\pm}}(q)=\nabla^2 \Psi^n_{\lambda_{\pm}}(q_0)(q-q_0)+o(\|q-q_0 \|_{\h^1_{2\pi}}).$ 
For  $\epsilon > 0$ sufficiently small, the homotopy 
$$H_{\pm} : (D(q_0,\epsilon) \times [0,1],\partial D(q_0,\epsilon) \times [0,1]) \to (\bH, \bH \setminus \{0\})$$ defined by 
$$H_{\pm}(q,\sigma)=\nabla^2 \Psi^n_{\lambda_{\pm}}(q_0)(q-q_0)+ \sigma  o(\|q-q_0 \|_{\h^1_{2\pi}})$$
 is well defined  $\sone$-equivariant gradient  homotopy, where $D(q_0,\epsilon)=\{q \in \bH : \|q-q_0\|_{\h^1_{2\pi}} \leq \epsilon\}$. 
Note that the $\sone$-invariant  potential $\Pi^n_{\lambda_{\pm}} : \bH^n \to \bR$ of $H_{\pm}(q,0)=\nabla^2 \Psi^n_{\lambda_{\pm}}(q_0)(q-q_0)$ is defined by $\ds \Pi^n_{\lambda_{\pm}}(q)=\frac{1}{2} \langle \nabla^2 \Psi^n_{\lambda_{\pm}}(q_0)(q-q_0),q-q_0\rangle_{\h^1_{2\pi}}$ i.e. $\nabla \Pi^n_{\lambda_{\pm}} : \bH^n \to \bH^n$ is a self-adjoint $S^1$-equivariant linear map, see formulas \eqref{fun} and  \eqref{fourier}. Since $\{q_0\}$ is an isolated zero along the homotopy $H(\cdot,\sigma),$ it follows that $\ci_{\sone}(\{q_0\},-\nabla \Psi^n_{\lambda_{\pm}})=\ci_{\sone}(\{q_0\},-\nabla \Pi^n_{\lambda_{\pm}}).$      In this way we significantly simplified computations of the $S^1$-equivariant Conley index $\ci_{\sone}(\{q_0\},-\nabla \Psi^n_{\lambda_{\pm}}).$ 

 To complete the proof it is enough to show that  there is $n_0 \in \bN$ such that for $n \geq n_0$  
 $$\ci_{\sone}(\{q_0\},-\nabla \Pi^{n}_{\lambda_{\pm}}) = \ci_{\sone}(\{q_0\},-\nabla \Pi^{n_0}_{\lambda_{\pm}})$$
  and that 
  $$\chi_{\sone}(\ci_{\sone}(\{q_0\},-\nabla \Pi^{n_0}_{\lambda_-})) \neq  \chi_{\sone}(\ci_{\sone}(\{q_0\},-\nabla \Pi^{n_0}_{\lambda_+})).$$
In the rest of the proof we show that these conditions are fulfilled.

Since  $q_0 \in \bH^1_{2\pi}$ is a constant function, $G(q_0)=\Gamma(q_0) \subset \mathbb{H}_0=\mathbb{R}^n \subset \bH^1_{2\pi}$ and 
$\ds 
\bH=T_{q_0}^{\perp}G(q_0)=\overline{T_{q_0}^{\perp} \Gamma(q_0) \oplus \bigoplus_{k=1}^{\infty} \bH_k}.
$

By assumption, for every $j=1,\ldots,j_0-1$ choose $k_j\in\bN$ such that $k_j^2<(\beta_j/\beta_{j_0})^2<(k_j+1)^2$ and note that $k_1 \geq k_2 \geq \ldots \geq k_{j_0-2} \geq k_{j_0-1}.$ 
Taking into account that $\frac{(k_1+1)^2}{\beta_1^2} > \frac{1}{\beta_{j_0}^2},$  $\lambda_+=\frac{1+\epsilon}{\beta_{j_0}}$ and that $\epsilon$  is arbitrarily small for fixed $n_0\geq k_1+1$ and $j=1,\ldots,m$ we obtain
\beq \label{szac}
n_0^2-\lambda_{\pm}^2\beta_j^2\geq n_0^2-\lambda_{\pm}^2\beta_1^2\geq n_0^2-\lambda_+^2\beta_1^2\geq \beta_1^2 \left(\frac{(k_1+1)^2}{\beta_1^2}-\lambda_{+}^2\right)>0.
\eeq

From the above formula and  equation (\ref{fourier}) it follows that for any $n \geq n_0$ the following equality holds 
$
\morse(-\nabla^2 {\Pi^n_{\lambda_{\pm}}})=\morse(-\nabla^2 {\Pi^{n_0}_{\lambda_{\pm}}}),
$ 
where $\morse(\cdot)$ is the Morse index.  Hence for any $n \geq n_0$ we obtain 
$\ci_{\sone}(\{q_0\},-\nabla {\Pi^n_{\lambda_{\pm}}})= \ci_{\sone}(\{q_0\},-\nabla {\Pi^{n_0}_{\lambda_{\pm}}}).$

%Therefore  the  $\sone$-equivariant Conley index of the isolated   invariant set $\{q_0\}$ under the vector field $-\nabla \Pi_{\lambda_{\pm}}$ is the $\sone$-homotopy type of a spectrum $(E_{n,\pm})_{n=n_0}^{\infty}$ (see \cite{[IZYDOREK]}), where $E_{n,\pm}$ is the same   pointed topological $\sone$-space for every $n\geq n_0.$ Consequently, the study of a change of the $\sone$-equivariant Conley index of $\{q_0\}$ under the  vector field $-\nabla\Pi_{\lambda_{\pm}}$ is equivalent to the study of finite-dimensional $\sone$-equivariant Conley index  $E_{n_0,\pm}=\ci_{\sone}(\{q_0\},-\nabla \Pi_{_{\lambda_{\pm}}|\h^{n_0}}).$

Since $\epsilon>0$ is arbitrarily small and $\lambda_{\pm}=\frac{1 \pm \epsilon}{\beta_{j_0}}$, for  $\lambda \in [\lambda_-,\lambda_+]$ we have 

\beq 
\label{prop1}
\text{if }  k>1  \text{ and } j=1,\ldots, m \quad  \text{ then }  \quad k^2-\lambda^2 \beta_j^2 \neq 0,
\eeq
\beq 
\label{prop2}
 \text{if }  k=1  \text{ and } j \neq j_0  \quad \text{ then } \quad k^2-\lambda^2 \beta_j^2 \neq 0,
\eeq
\beq 
\label{prop3}
\text{if }  k=1,  j=j_0   \text{ and } k^2-\lambda^2 \beta_j^2 = 0 \quad \text{ then } \quad \lambda=\frac{1}{\beta_{j_0}}.
\eeq
Moreover, it is clear that 
\beq \label{prop4}(1-\lambda^2_- \beta_{j_0}^2)(1-\lambda^2_+ \beta_{j_0}^2)=-\epsilon^2(4-\epsilon^2)<0.
\eeq
Applying formulas \eqref{prop2}, \eqref{prop4} we obtain that 
\beq \label{prop5} \morse\left(-\nabla^2 {\Pi^1_{\lambda_-}}\right) \neq \morse\left(-\nabla^2 {\Pi^1_{\lambda_+}}\right),
\eeq where $\morse(\cdot)$ is  the Morse index.

Taking into account formulas \eqref{fourier}, \eqref{prop1}   and the spectral decomposition of $\h^{n_0}$ given by the isomorphisms $-\nabla^2 \Pi^{n_0}_{{\lambda_{\pm}}}(q_0)$ we obtain 
$$\bH^{n_0}=\bH_1 \oplus \left(T_{q_0}^{\perp} \Gamma(q_0)\oplus \bigoplus_{k=2}^{n_0} \bH_{k} \right)= \bH_1  \oplus \bW = \bH_1 \oplus \bW^- \oplus \bW^+.$$
By formulas \eqref{fourier}, \eqref{prop2}, \eqref{prop3}   and the spectral decomposition of $\bH^{n_0}$ given by the isomorphisms $-\nabla^2 \Pi^{n_0}_{{\lambda_{+}}}(q_0)$ we get 
$$\bH^{n_0}=\bH_1 \oplus \left(T_{q_0}^{\perp} \Gamma(q_0)\oplus \bigoplus_{k=2}^{n_0} \bH_{k} \right)= \left(\bH_{1,+}^- \oplus \bH_{1,+}^+\right) \oplus \left(\bW^- \oplus \bW^+\right).$$
Moreover the spectral decomposition of $\bH^{n_0}$ given by the isomorphisms $-\nabla^2 \Pi^{n_0}_{\lambda_-}(q_0)$ is of the form 
 $$\bH^{n_0}=\bH_1 \oplus \left(T_{q_0}^{\perp} \Gamma(q_0)\oplus \bigoplus_{k=2}^{n_0} \bH_{k} \right)= \left(\bH_{1,-}^- \oplus \bH_{1,-}^+\right) \oplus \left(\bW^- \oplus \bW^+\right).$$

\nt It follows that $\ds \ci_{\sone} (\{q_0\},-\nabla {\Pi^{n_0}_{\lambda_{\pm}}})= S^{\bH_{1,\pm}^+} \wedge S^{\bW^+}.$ Hence we obtain 

\beq \label{prop6} \chi_{\sone} \left( \ci_{\sone} (\{q_0\},-\nabla {\Pi^{n_0}_{\lambda_{\pm}}})\right)= \chi_{\sone}\left(S^{\bH_{1,{\pm^+}}}\right) \star  \chi_{\sone}\left( S^{\bW^+}\right) \in U(\sone).
\eeq
We claim that $\ds \chi_{\sone}\left(\ci_{\sone} (\{q_0\},-\nabla {\Pi^{n_0}_{\lambda_-}})\right) \neq \chi_{\sone}\left(\ci_{\sone} (\{q_0\},-\nabla {\Pi^{n_0}_{\lambda_{+}}})\right).$
Suppose contrary to our claim that $\ds \chi_{\sone}\left(\ci_{\sone} (\{q_0\},-\nabla {\Pi^{n_0}_{\lambda_-}})\right) = \chi_{\sone}\left(\ci_{\sone} (\{q_0\},-\nabla {\Pi^{n_0}_{\lambda_{+}}})\right).$ By Remark \ref{invert}  the element $\chi_{\sone}(S^{\bW^+})$ is invertible in the Euler ring  $U(\sone).$ Therefore taking into account formula \eqref{prop6} we obtain   $\chi_{\sone}\left(S^{\bH_{1,-}^+}\right)= \chi_{\sone}\left(S^{\bH_{1,+}^+}\right).$ By formula \eqref{prop5} we have $r_-:=\dim \bH^+_{1,-} \slash 2 \neq \dim \bH^+_{1,+}\slash 2 =: r_+.$

Since $\bH_1= \lin \{e_i \cos t, e_i \sin t: i=1,\ldots,n\}$  and the action of the group $S^1$ on $\bH_1$ is given by shift in time, the spaces $\bH^+_{1,\pm}$  are representations of the group $\sone$ such that $\bH^+_{1,\pm} \approx_{\sone} \bR[r_{\pm},1].$ Hence by formula \eqref{ches} we have 
$$\chi_{\sone}\left(S^{\bH^+_{1,\pm}}\right)= \chi_{\sone}\left(S^{\bR[r_{\pm},1]}\right)=    \bI - r_{\pm} \chi_{\sone}\left({\sone \slash \bZ_1}^+\right) \in U(\sone).$$
It follows that $r_-=r_+,$ a contradiction.

We have just proved that $\ds \chi_{\sone}\left(\ci_{\sone} (\{q_0\},-\nabla {\Pi^{n_0}_{\lambda_-}})\right) \neq \chi_{\sone}\left(\ci_{\sone} (\{q_0\},-\nabla {\Pi^{n_0}_{\lambda_+}})\right).$
Hence for $n \geq n_0$ $$E_{n,-}=\ci_{\sone} (\{q_0\},-\nabla {\Psi^n_{\lambda_-}})=\ci_{\sone} (\{q_0\},-\nabla {\Pi^{n_0}_{\lambda_-}}) \not \approx_{\sone}$$
$$\ci_{\sone} (\{q_0\},-\nabla {\Pi^{n_0}_{\lambda_+}}) =\ci_{\sone} (\{q_0\},-\nabla {\Psi^n_{\lambda_+}})=E_{n,+}.$$ In other words for any $n \geq n_0$ the spaces $E_{n,-}$ and $E_{n,+}$ are not $S^1$-homotopically equivalent. Consequently the $S^1$-homotopy types of spectra $(E_{n,-})_{n=n_0}^{\infty}$ and $(E_{n,+})_{n=n_0}^{\infty}$ (see \cite{[IZYDOREK]}) are different. Hence $\ci_{\sone} (\{q_0\},-\nabla \Psi_{\lambda_-}) \neq \ci_{\sone} (\{q_0\},-\nabla \Psi_{\lambda_+}),$ which completes the proof.
\end{proof}

\subsection{Proof of Theorem \ref{main-theo}}
\begin{proof}
It is well known for people in the field that a change of the $G=(\Gamma \times \sone)$-equivariant Conley index  $\cig(G(q_0),-\nabla \Phi(\cdot,\lambda)) $ along the family $\cT$ implies the existence of a local bifurcation of  solutions of equation  \eqref{gradf} from this family. Therefore in order to obtain a bifurcation of solutions of equation  \eqref{gradf} from the orbit $G(q_0) \times \{1 \slash \beta_{j_0}\} \subset \bH^1_{2 \pi} \times (0,+\infty)$   it is enough to show that   $\cig(G(q_0),-\nabla \Phi(\cdot,\lambda_-)) \neq \cig(G(q_0),-\nabla \Phi(\cdot,\lambda_+)),$ where  $\lambda_{\pm}=\frac{1 \pm \epsilon} {\beta_{j_0}}$ and $\epsilon>0$ is sufficiently small. This inequality will be a consequence of lemmas \ref{gh}, \ref{cis} and corollary \ref{hgg}. 

\nt Recall that $\ds \bH^n= T_{q_0}^{\perp} \Gamma(q_0) \oplus \bigoplus_{k=1}^{n} \bH_k.$ 
By Lemma \ref{cis} we obtain $n_0 \in \bN$ such that  for every $n \geq n_0$
\beq \label{ineq} \chi_{\sone}(\ci_{\sone}(\{q_0\},-\nabla \Psi^n_{\lambda_-})) = \chi_{\sone}(\ci_{\sone}(\{q_0\},-\nabla \Psi^{n_0}_{\lambda_-})) \neq $$ $$\chi_{\sone}(\ci_{\sone}(\{q_0\},-\nabla \Psi^{n_0}_{\lambda_+}))= \chi_{\sone}(\ci_{\sone}(\{q_0\},-\nabla \Psi^n_{\lambda_+})).
\eeq
By inequality \eqref{szac}  for any $n \geq n_0$ the following equality holds 
$$
\morse(\nabla^2 \Psi^n_{\lambda_{\pm}}(q_0))=\morse(-\nabla^2 {\Pi^n_{\lambda_{\pm}}})=\morse(-\nabla^2 {\Pi^{n_0}_{\lambda_{\pm}}}),
$$
where $\morse(\cdot)$ is the Morse index.  

\nt It follows that for any $n \geq n_0$ we obtain 
$\cig(G(q_0),-\nabla \Phi^n(\cdot, \lambda_{\pm}))= \cig(G(q_0),-\nabla \Phi^{n_0}(\cdot, \lambda_{\pm})),$ where $ \Phi^n(\cdot, \lambda_{\pm}) = \Phi(\cdot, \lambda_{\pm})_{\mid \oplus_{k=0}^n \bH_k} :\ds \bigoplus_{k=0}^n \bH_k \to \bR$ (see \cite{[IZYDOREK]} for more details).

 Since the isotropy group $\Gamma_{q_0}$ is trivial, the isotropy group $G_{q_0} \in \subg$ of $q_0 \in \bH^1_{2 \pi}$ equals $\{e\} \times \sone.$ By Lemma \ref{gh} the pair $(G,G_{q_0})=(\Gamma \times \sone,\{e\} \times \sone)$ is admissible.  Hence by Theorem \ref{cio},  condition \eqref{ineq} and Corollary \ref{hsmadmco} we obtain that  for $n \geq n_0$ 
 $$\chig(\cig(G(q_0),-\nabla \Phi^n(\cdot, \lambda_-))) = \chig(G \wedge_{\sone} \ci_{\sone}(\{q_0\},-\nabla \Psi^n_{\lambda_-})) \neq $$ $$\chig(G \wedge_{\sone} \ci_{\sone}(\{q_0\},-\nabla \Psi^n_{\lambda_+})) =\chig(\cig(G(q_0),-\nabla \Phi^n(\cdot, \lambda_{+})))$$ 
 and consequently  $\cig(G(q_0),-\nabla \Phi^n(\cdot, \lambda_-)) \neq \cig(G(q_0),-\nabla \Phi^n(\cdot, \lambda_+)),$ for any $n \geq n_0.$   

\nt Therefore  the  $G$-equivariant Conley index of the isolated   invariant set $G(q_0)$ under the vector field $-\nabla \Phi(\cdot, \lambda_{\pm})$ is the $G$-homotopy type of a spectrum $(\cE_{n,\pm})_{n=n_0}^{\infty}$ (see \cite{[IZYDOREK]}), where $\cE_{n,\pm}=\cig(G(q_0),-\nabla \Phi^n(\cdot,\lambda_{\pm}))$ is the same   pointed topological $G$-space for every $n\geq n_0.$ 

\nt Summing up, we have just proved that $\cig(G(q_0),-\nabla \Phi(\cdot,\lambda_-)) \neq \cig(G(q_0),-\nabla \Phi(\cdot,\lambda_+)),$ which completes the proof.
\end{proof}

\br 
Modifying slightly the proof of Theorem  \ref{main-theo}   one can show that a connected set of non-stationary periodic solutions of system \eqref{newsys} emanate from the orbit $\Gamma(q_0).$ In the proof of Theorem \ref{main-theo} we have shown that  the  $G$-equivariant Conley index of the isolated   invariant set $G(q_0) = \Gamma(q_0) \subset \bH_0=\bR^n \subset \bH^1_{2\pi}$ under the vector field $-\nabla \Phi(\cdot, \lambda_{\pm})$ i.e.  $\cig(G(q_0),-\nabla \Phi(\cdot,\lambda_{\pm}))$ (see \cite{[IZYDOREK]}), is the $G$-homotopy type of spectrum $(\cE_{n,\pm})_{n=n_0}^{\infty}$  where $\cE_{n,\pm}=\cig(G(q_0),-\nabla \Phi^n(\cdot, \lambda_{\pm}))$ for every $n\geq n_0$ i.e. this spectrum is constant. 

Let $\Upsilon_G(\cdot)$ be the $G$-equivariant Euler characteristic for $G$-homotopy types of $G$-equivariant spectra defined in   \cite{[GORY1]}. Since the operator  $\nabla \Phi(\cdot,\lambda_{\pm})$ is of the form compact perturbation of the identity, directly from the definition of $\Upsilon_G(\cdot)$  it  follows that $\Upsilon_G(\cig(G(q_0),-\nabla \Phi(\cdot,\lambda_{\pm})))=$ $\chig(\cig(G(q_0),-\nabla \Phi^{n_0}(\cdot,\lambda_{\pm}))) \in U(G).$  

Let  $\cO \subset \bH^1_{2\pi}$ be an open bounded and $G$-invariant subset such that $\nabla \Phi(\cdot,\lambda_{\pm})^{-1}(0) \cap \cO = G(q_0).$ It was shown  in \cite{[GORY1]}  that  $\Upsilon_G(\cig(G(q_0),-\nabla \Phi(\cdot,\lambda_{\pm})))= \nabla_G\textrm{-}\mathrm{deg}(\nabla \Phi(\cdot,\lambda_{\pm}), \cO) \in U(G)$ is the degree for $G$-equivariant gradient maps defined in \cite{[GORY]}. From the proof of Theorem \ref{main-theo} it follows that $\chig(\cig(G(q_0),-\nabla \Phi^{n_0}(\cdot, \lambda_{-}))) \neq \chig(\cig(G(q_0),-\nabla \Phi^{n_0}(\cdot, \lambda_{+}))).$ 

Hence $\nabla_G\textrm{-}\mathrm{deg}(\nabla \Phi(\cdot,\lambda_{-}), \cO) \neq \nabla_G\textrm{-}\mathrm{deg}(\nabla \Phi(\cdot,\lambda_{+}), \cO).$ It is known that a change of this degree implies bifurcation of a connected set of solutions of equation $\nabla \Phi(q,\lambda)=0$ from the orbit $\Gamma(q_0) \times \{1 \slash \beta_{j_0}\} \subset \bH^1_{2 \pi} \times (0,+\infty).$

Finally we would like to underline that basic material on degree theories for equivariant maps can be found in \cite{[BKRS],[BAKRST],[RYBICKI]}.
\er

%%%%%%%%%%%%%%%%%%%%%%%%%%%%%%%%%%%%%%%%%%%%%%%%%%%%%%%

\section{Applications}\label{applications}
In order to show the strength of our main result \ref{main-theo} in this section we apply it to a couple of special class of galactic potentials. We point out that our goal here is not the analysis of specific galaxies and their dynamics, we just want take  a kind of generic galactic potentials to show how to find periodic orbits on them. Must of the work on galactic potentials is numeric, we will show here, just in a couple of simple cases, the way to obtain periodic orbits in an analytic way. The target is that it could be used as a started point for people working in the field.

Since many galactic potentials defined in the plane are of the form $ U(x^2,y^2)$ (see for instance \cite{[ALF],[PUC]}), or more generally must of them are considered as a perturbation of an harmonic oscillator (see \cite{[CAR]} for more details) having the form $$ U(x,y)= \omega^2 (x^2+y^2) + \varepsilon V(||(x,y)||^2)$$ and since we have studied symmetric potentials along this paper, we will apply the symmetric Liapunov center theorem to  potentials of the form $U(||(x,y)||^2)$, which can be seen as a special class of galactic potentials.

 We believe that generalizations of important results as the Liapunov center theorem for symmetric potentials could be useful in some applications (see for instance the generalization of the famous Weierstrass model for homogeneous potentials \cite{[BK]})
 
\bex
Assume that a galaxy on the plane is moving under the influence of the potential $U(x) = -2\|x\|^4 + \frac{5}{3} \|x\|^6 - \frac{1}{4} \|x\|^8.$ 
Consider the polynomial $\vp : \bR \to \bR$  defined by $\vp (t) = -2t^2 + \frac{5}{3}t^3 - \frac{1}{4}t^4.$ We observe that  $\vp'(t)=-t (t-1) (t-4)$ and $U(x)=\vp(\|x\|^2).$

Since the gradient $\nabla U : \bR^2 \to \bR^2$ is given by $\nabla U(x)=2 \vp'(\|x\|^2)x,$ we obtain

$$(\nabla U)^{-1}(0)=\{(0,0)\} \cup \{x \in \bR^n : \|x\|=1\} \cup \{x \in \bR^n : \|x\|=2\}=$$ $$= S_0 \cup S_1 \cup S_2=\sotwo (0,0) \cup \sotwo (1,0) \cup \sotwo (2,0).$$

Taking into account that  $\nabla_{x_i}U(x) = 2 \vp'(\|x\|^2)x_i$ we compute the Hessian 

$$\nabla^2  _{x_jx_i}U(x) = \left\{\ba{cc}  4 \vp''(\|x\|^2)x_i x_j & \textrm{ for } i \not = j, \\  & \\ 2 \vp'(\|x\|^2) +  4 \vp''(\|x\|^2)x^2_i & \textrm{ for } i  = j. \ea.  \right. $$

Since $\vp''(t)=-3t^2+10t-4$  we obtain:
\begin{itemize}
\item for $(0,0) \in S_0, \quad \nabla^2 U((0,0))=\left[\ba{cc}  0  & 0 \\ 0 &  0 \ea  \right].$
\item for $(1,0) \in S_1,  \quad \nabla^2 U((1,0))=\left[\ba{cc}  12  & 0 \\ 0 &  0 \ea  \right].$
\item for $(2,0) \in S_2$, \quad $\nabla^2 U((2,0))=\left[\ba{cc}  -192  & 0 \\ 0 &  0 \ea  \right].$
\end{itemize}

It is easy verify that the hypothesis  of  Theorem \ref{main-theo} are fulfilled at orbit $S_1=\sotwo (1,0)$ but not at orbit $S_2=\sotwo (2,0)$. Therefore in any neighborhood of $S_1=\sotwo (1,0)$ there exists at least one periodic orbit.
\eex

\bex
In general we can think in the interaction of several galaxies moving under the 
influence of a $\sotwo$-invariant  potential, in this way we define  
 polynomials $\cU=\cU(t_1,\ldots,t_m) : \bR^m \to \bR$ and $\cU_0 : \bR^{2m} \to \bR$. Now we define an $\sotwo$-invariant  potential   $U : \bR^{2m} \to \bR$  by 
 \begin{eqnarray*}
 U(x) &=& \frac{\omega^2}{2} \|x\|^2 + \frac{\varepsilon}{2} \cU_0 (x_1^2,x_2^2,\ldots,x_{2m-1}^2, x_{2m}^2) \\
 &=& \frac{\omega^2}{2} \|x\|^2 + \frac{\varepsilon}{2} \cU (x_1^2+x_2^2,\ldots,x_{2m-1}^2 + x_{2m}^2).  \end{eqnarray*} 

We have used the auxiliary polynomial $\cU_0$ to clarify that it depends only on the squares of the variables to note the similarity with the known polynomial galactic potentials (see \cite{[CAR]}). 
Note that for $i=1,\ldots,m$ we have
%$$\nabla_{x_{2i-1}}U(x)= x_{2i-1}(\omega^2 + \varepsilon \cU'_{t_i}(x_1^2+x_2^2,\ldots,x_{2m-1}^2 + x_{2m}^2)),$$
%$$\nabla_{x_{2i}}U(x)= x_{2i}(\omega^2 + \varepsilon \cU'_{t_i}(x_1^2+x_2^2,\ldots,x_{2m-1}^2 + x_{2m}^2)).$$
\[
\frac{\partial U}{\partial x_{2i-1}}(x)= x_{2i-1}(\omega^2 + \varepsilon \frac{\partial\cU}{\partial t_i}(x_1^2+x_2^2,\ldots,x_{2m-1}^2 + x_{2m}^2)),
\]
\[
\frac{\partial U}{\partial x_{2i}}(x)= x_{2i}(\omega^2 + \varepsilon \frac{\partial\cU}{\partial t_i}(x_1^2+x_2^2,\ldots,x_{2m-1}^2 + x_{2m}^2)).
\]
Taking into consideration the above we obtain that $\nabla U(x)=0$ iff
\[
\forall_{i=1,\ldots,m} \left(\omega^2 + \varepsilon \frac{\partial\cU}{\partial t_i}(x_1^2+x_2^2,\ldots,x_{2m-1}^2 + x_{2m}^2)=0 \quad \vee \quad x_{2i-1}=x_{2i}=0\right).
\]
Moreover, putting $q(x)=(x_1^2+x_2^2,\ldots,x_{2m-1}^2 + x_{2m}^2)$, we have
\begin{eqnarray*}
\frac{\partial^2 U}{\partial x_{k}^2}(x) &=& \left(\omega^2 + \varepsilon \frac{\partial\cU}{\partial t_i}\left(q(x)\right)\right)+x_{k}^2\left(2\varepsilon\frac{\partial^2\cU}{\partial^2 t_i}\left(q(x)\right)\right), \\
\frac{\partial^2 U}{\partial x_{n}\partial x_{k}}(x) &=& x_{n}x_{k}\left(2\varepsilon \frac{\partial^2\cU}{\partial t_i\partial t_j}\left(q(x)\right)\right),
\end{eqnarray*}
where $k\in\{2i-1,2i\}$, $n\in\{2j-1,2j\}$ and  $k\neq n$.

We apply the above to a concrete simple case. Let be $m=2$ and $\ds \cU(t_1,t_2)=- \frac{1}{2} t_1^2 + \frac{1}{2} t_1^2 t_2^4$. In this case $\nabla U(x)=0$ iff

$$ \omega^2 - \varepsilon (x_1^2+x_2^2)  + \epsilon (x_1^2+x_2^2)(x_3^2+x_4^2)^4=0 \text{ or } x_1=x_2=0$$
and
$$ \omega^2    + 2\varepsilon (x_1^2+x_2^2)^2(x_3^2+x_4^2)^3=0 \text{ or } x_3=x_4=0.$$
After straightforward computations we get $\nabla U(x)=0 \text{ iff }$ $\ds x_1=x_2=x_3=x_4=0 \text{ or }x_1^2+x_2^2=\frac{\omega^2}{\varepsilon},x_3=x_4=0$ i.e. 
\begin{eqnarray*}
\nabla (U)^{-1}(0) &=& \{(0,0,0,0)\} \cup \sotwo (\frac{\omega}{\sqrt{\varepsilon}},0,0,0) \\ &=& \{(0,0,0,0)\} \cup \sotwo (\frac{\omega}{\sqrt{\varepsilon}},0) \times \{(0,0)\}.\end{eqnarray*}
Computing the Hessian $\nabla^2 U(x)$ at $x'=(0,0,0,0)$ and $\ds x''=(\frac{\omega}{\sqrt{\varepsilon}},0,0,0)$ we obtain

$$\nabla^2 U(x')=\left[\ba{cccc} \omega^2 & 0 & 0 & 0 \\ 0 & \omega^2 & 0 & 0 \\ 0 & 0 & \omega^2 & 0 \\
0 & 0 & 0 & \omega^2\ea\right] 
\quad \text{ and } \quad
\nabla^2 U(x'')=\left[\ba{cccc} -2\omega^2 & 0 & 0 & 0 \\ 0 & 0 & 0 & 0 \\ 0 & 0 & \omega^2 & 0 \\
0 & 0 & 0 & \omega^2\ea\right].$$ 

\bigskip

We verify that $x'$ satisfies assumptions of the classical Liapunov center theorem, which means that the origin bifurcate into a family of periodic orbits. We also observe that we can not apply this theorem to the point $x''$ because the Hessian at this point is degenerate, nevertheless we can apply our main Theorem \ref{main-theo}, the symmetric Liapunov center theorem to the orbit $\sotwo(x'')$, getting periodic solutions in any neighborhood of that orbit, i.e. 
we have a local bifurcation of $\sotwo(x'')$ into periodic orbits.
\eex

 \end{document}